\documentclass[leqno]{article}
\usepackage{latexsym,amsmath,amsfonts,mathrsfs,amssymb,amsthm}
\usepackage[usenames]{color}
\usepackage{graphicx}
\usepackage{setspace}
\makeatletter
\renewcommand\normalsize{%
   \@setfontsize\normalsize\@xpt\@xiipt
   \abovedisplayskip 3\p@ \@plus2\p@ \@minus5\p@
   \abovedisplayshortskip \z@ \@plus3\p@
   \belowdisplayshortskip 6\p@ \@plus3\p@ \@minus3\p@
   \belowdisplayskip \abovedisplayskip
   \let\@listi\@listI}
\makeatother

\def\R{\mathbb R}
\def\N{\mathbb N}
\def\C{\mathbb C}

\def\a{\alpha}
\def\e{\epsilon}
\def\d{\delta}
\def\Y{\mathbb Y}
\def\T{\mathbb T}
\def\P{\mathbb P}
\def\be{\begin{equation}}
\def\ee{\end{equation}}
\def\bs{\backslash}
\def\qed{\hfill$\Box$\bigskip}
\def\nd{\noindent Proof. }
\numberwithin{equation}{section}
\newtheorem{lem}[equation]{Lemma}
\newtheorem{pro}[equation]{Proposition}
\newtheorem{defn}[equation]{Definition}
\newtheorem{thm}[equation]{Theorem}

\newtheorem{rem}[equation]{Remark}
\newtheorem*{thmA}{Theorem 7.1}
\newtheorem*{thmB}{Theorem 7.2}
\begin{document}
\bigskip

\centerline{\Large \textbf{Global regularity for minimal sets near a union of two planes}}

\bigskip

\centerline{\large Xiangyu Liang}

\bigskip

\centerline{ D\'epartement de Math\'ematiques, B\^at 425, Universit\'e Paris-Sud 11, 91405 Orsay Cedex}
\centerline{xiangyuliang@gmail.com}

\vskip 1cm

\centerline {\large\textbf{Abstract.}} 

\bigskip

We discuss the global regularity of 2 dimensional minimal sets that are near a union of two planes, and prove that every global minimal set in $\R^4$ that looks like a union of two almost orthogonal planes at infinity is a cone. The main point is to use the topological properties of a minimal set at a large scale to control its behavior at smaller scales. 
%

\bigskip


\bigskip


\section{Introduction}

This paper deals with the local (resp. global) regularity of two-dimensional minimal sets in $\R^4$ that looks like the union of two almost orthogonal planes locally (resp. at infinity). The motivation is that we want to decide whether all global minimal sets in $\R^n$ are cones.

This Bernstein type of problem is of typical interest for all kinds of minimizing problems in geometric measure theory and calculus of variations. It is natural to ask how does a global minimizer look like, as soon as we know already the local regularity for minimizers. Well known examples are the global regularity for complete 2-dimensional minimal surfaces in $\R^3$, area or size minimizing currents in $\R^n$, or global minimizers for the Mumford-Shah functional. Some of them admit very good descriptions. See \cite{Be, Mo89, Mo86,DMS} for further information.

Here our notion of minimality is defined in the setting of sets. Roughly speaking, we say that a set $E$ is minimal when there is no deformation $F=\varphi(E)$, where $\varphi$ is Lipschitz and $\varphi(x)-x$ is compactly supported, for which the Hausdorff measure $H^2(F)$ is smaller than $H^2(E)$. More precisely, 

\begin{defn}[Almgren competitor (Al competitor for short)] Let $E$ be a closed set in an open subset $U$ of $\R^n$ and $d\le n-1$ be an integer. An Almgren competitor for $E$ is a closed set $F\subset U$ that can be written as $F=\varphi_1(E)$, where $\varphi_t:U\to U$ is a family of continuous mappings such that 
\be \varphi_0(x)=x\mbox{ for }x\in U;\ee
\be\mbox{ the mapping }(t,x)\to\varphi_t(x)\mbox{ of }[0,1]\times U\mbox{ to }U\mbox{ is continuous;}\ee
\be\varphi_1\mbox{ is Lipschitz,}\ee
  and if we set $W_t=\{x\in U\ ;\ \varphi_t(x)\ne x\}$ and $\widehat W=\bigcup_{t\in[0.1]}[W_t\cup\varphi_t(W_t)]$,
then
\be \widehat W\mbox{ is relatively compact in }U.\ee
Such a $\varphi_1$ is called a deformation in $U$, and $F$ is also called a deformation of $E$ in $U$.
\end{defn}

\begin{defn}[(Almgren) minimal sets]
Let $0<d<n$ be integers, $U$ an open set of $\R^n$. A closed set $E$ in $U$ is said to be (Almgren) minimal of dimension $d$ in $U$ if 
\begin{spacing}{1}\be H^d(E\cap B)<\infty\mbox{ for every compact ball }B\subset U,\ee
and
\be H^d(E\bs F)\le H^d(F\bs E)\ee\end{spacing}
for all Al competitors $F$ for $E$.
\end{defn}

This notion was introduced by Almgren to modernize Plateau's problem, which aims at understanding physical objects, such as soap films, that minimize the area while spanning a given boundary.  The study of regularity and existence for these sets is one of the canonical interests in geometric measure theory.
%
%

\medskip

Our goal is to show that every minimal set in $\R^n$ is a cone. The general idea is the following.

Let $E$ be a $d-$dimensional reduced Almgren minimal set in $\R^n$. Reduced means that there is no unnecessary points. More precisely, we say that $E$ is reduced when
\be H^d(E\cap B(x,r))>0\mbox{ for } x\in E\mbox{ and }r>0.\ee

Recall that the definition of minimal sets is invariant modulo sets of measure zero, and it is not hard to see that for each Almgren (resp. topological) minimal set $E$, its closed support $E^*$ (the reduced  set $E^*\subset E$ with $H^2(E\bs E^*)=0$) is a reduced Almgren (resp. topological) minimal set. Hence we can restrict ourselves to discussing only reduced minimal sets.

Now fix any $x\in E$, and set
\be \theta_x(r)=r^{-d}H^d(E\cap B(x,r)).\ee
This density function $\theta_x$ is nondecreasing for $r\in ]0,\infty[$ (cf.\cite{DJT} Proposition 5.16). In particular the two values
\be \theta(x)=\lim_{t\to 0^+}\theta_x(t)\mbox{ and }\theta_\infty(x)=\lim_{t\to \infty}\theta_x(t)\ee
exist, and are called density of $E$ at $x$, and density of $E$ at infinity respectively. It is easy to see that $\theta_\infty(x)$ does not depend on $x$, hence we shall denote it by $\theta_\infty$.

Theorem 6.2 of \cite{DJT} says that if $E$ is a minimal set, $x\in E$, and $\theta_x(r)$ is a constant function of $r$, then $E$ is a minimal cone centered on $x$. Thus by the monotonicity of the density functions $\theta_x(r)$ for any $x\in E$, if we can find a point $x\in E$ such that $\theta(x)=\theta_\infty$, then $E$ is a cone and we are done.

On the other hand, the possible values for $\theta(x)$ and $\theta_\infty$ for any $E$ and $x\in E$ are not arbitrary. 
By Proposition 7.31 of \cite{DJT}, for each $x$, $\theta(x)$ is equal to the density at the origin of a $d-$dimensional Al-minimal cone in $\R^n$. An argument around (18.33) of \cite{DJT}, which is similar to the proof of Proposition 7.31 of \cite{DJT}, gives that $\theta(x)$ is also equal to the density at the origin of a $d-$dimensional Al-minimal cone in $\R^n$. In other words, if we denote by $\Theta_{d,n}$ the set of all possible numbers that could be the density at the origin of a $d-$dimensional Almgren-minimal cone in $\R^n$, then $\theta_\infty\in\Theta_{d,n}$, and for any $x\in E$, $\theta(x)\in\Theta_{d,n}$.

Thus we restrict the range of $\theta_\infty$ and $\theta(x)$. Recall that the set $\Theta_{d,n}$ is possibly very small for any $d$ and $n$. For example, $\Theta_{2,3}$ contains only three values: 1 (the density of a plane), 1.5 (the density of a $\Y$ set, which is the union of three closed half planes with a common boundary $L$, and that meet along the line $L$ with $120^\circ$ angles), and $d_T$ (is the density of a $\T$ set, i.e., the cone over the 1-skeleton of a regular tetrahedron centered at 0). (See the figure below).

\centerline{\includegraphics[width=0.16\textwidth]{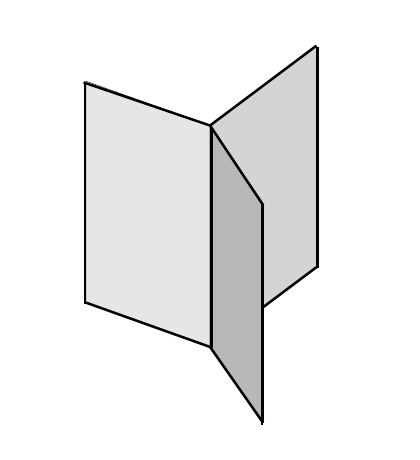} \hskip 1.5cm\includegraphics[width=0.2\textwidth]{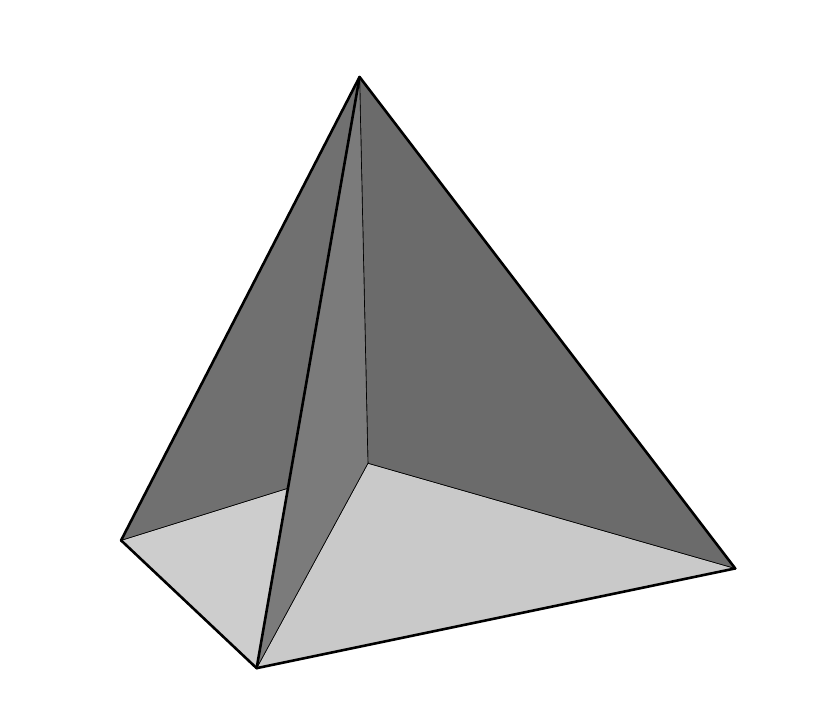}}
\hskip 4.1cm a $\Y$ set\hskip 3.7cm  a $\T$ set

Recall that the reason why $\theta_\infty$ has to lie in $\Theta_{d,n}$ is that, for any Al-minimal set $E$, all its blow-in limits have to be Al-minimal cones (cf. Argument around (18.33) of \cite{DJT}). A blow-in limit of $E$ is the limit of any converging (for the Hausdorff distance) subsequence of
\be E_r=r^{-1}E, r\to\infty .\ee

Hence the value of $\theta_\infty$ implies that at sufficiently large scales, $E$ looks like an Al-minimal cone of density $\theta_\infty$. 

This is the same reason why $\theta(x)\in\Theta_{d,n}$. Here we look at the behavior of $E_r$ when $r\to 0$, and the limit of any converging subsequence is called a blow-up limit (this might not be unique!). Such a limit is also an Al-minimal cone $C$ (cf. \cite{DJT} Proposition 7.31). This means, at some very small scales around each $x$, $E$ looks like some Al-minimal cone $C$ of density $\theta(x)$. In this case we call $x$ a $C$ type point of $E$.

After the discussion above, our problem will be solved if we can prove that every minimal cone $C$ satisfies the following property: 
\be\begin{split}&\mbox{There exists }\e=\e_C>0\mbox{, such that for every minimal set }E\mbox{, if }d_{0,1}(C,E)<\e
\mbox{, then }\\&\mbox{there exists }x\in E\cap B(0,1)\mbox{ whose density }\theta(x)\mbox{ is the same as that of }C\mbox{ at the origin.}\end{split}\ee

Here $d_{x,r}$ stands for the relative distance in the ball $B(x,r)$: for any closed sets $E$ and $F$,
\be d_{x,r}(E,F)=\frac1r\max\{\sup\{d(y,F):y\in E\cap B(x,r)\},\sup\{d(y,E):y\in F\cap B(x,r)\}\}.\ee
%
%
%

The discussion above uses only the values of densities at small scale and at infinity. A geometric intepretation is: there exists $x\in E\cap B(0,1)$ such that 
a blow-up limit $C_x$ of $E$ at $x$ admits the same density as $C$ at the origin. 
%
%
%

\medskip

So far we know that (1.13) is true for the planes and $\Y$ sets (see \cite{DJT} Proposition 16.24). We do not know any minimal cone that does not verify the property (1.13). But there are at least two minimal cones for which we do not know whether (1.13) holds, either: the $\T$ set, and the sets $Y\times Y\in \R^4$, whose minimality has recently been proved in \cite{YXY}. The topology of the set $Y\times Y$ is more complicated than that of $\T$ sets, and the situation of $\T$ sets is already tricky, see \cite{globalT} for more detail.

In this paper we prove the property (1.13) for the unions of two almost orthogonal planes. Recall that in \cite{2p}, we have proved the following
\begin{thm}[minimality of the union of two almost orthogonal planes, cf. \cite{2p} Thm 1.24]There exists $0<\theta_0<\frac\pi2$, such that if $P^1$ and $P^2$ are two planes in $\R^4 $ whose characteristic angles $(\alpha_1,\alpha_2)$ satisfy $\a_2\ge\a_1\ge \theta$, then their union $P^1\cup P^2$ is a minimal cone in $\R^4$.
\end{thm}

Here the characteristic angles describe the relative position between planes. Two planes $P^1$ and $P^2$ have characteristic angles $(\alpha_1,\alpha_2)$ with $\a_2\ge\a_1\ge \theta$ means that there exists an orthonormal basis $\{e_i\}_{1\le i\le 4}$ of $\R^4$ such that $P_\a^1$ is generated by $e_1$ and $e_2$, and $P_\a^2$ is generated by $\cos\a_1e_1+\sin\a_1 e_3$ and $\cos\a_2e_2+\sin\a_2 e_4$. Each pair of $\a=(\a_1,\a_2)$ with $\a_2\ge \a_1\ge\theta$ gives a minimal cone $P_\a=P^1\cup_\a P^2$, and the origin is called a singularity of type $\P_\a$ in the set $P_\a$. These gives a continuous family of minimal cones with the same density at the origin, any two of which are not $C^1$ equivalent to each other. But still, we give them a general name, that is, each singularity of type $\P_\a$ is a singular point of type $2\P$.

So let us state our main results. 

\begin{thm}\label{main}There is an angle $\theta_1\in[\theta_0,\frac\pi2)$, (where $\theta_0$ is the $\theta_0$ in Theorem 1.15), and $\lambda>0$, such that for any $\a=(\a_1,\a_2)$ with $\a_2\ge\a_1\ge\theta_1$, if $E$ is a 2-dimensional reduced Almgren minimal set in $U\subset\R^4$, $B(x,r)\subset U$, and there is a reduced minimal cone $P_\a$ of type $\P_\a$ centered at $x$ such that $d_{x,r}(E,P_\a)\le\lambda$, then $E\cap B(x,r/100)$ contains (at least) a $2\P$ type point.
\end{thm}

A direct corollary to this is the expected global regularity for minimal sets that look like a union of two plane at the infinity:

\begin{thmA}Let $\theta_1$ be as in Theorem \ref{main}. Then for any $\a=(\a_1,\a_2)$ with $\a_2\ge \a_1\ge\theta_1$, if $E$ is a 2-dimensional reduced Almgren minimal set in $\R^4$ such that one blow-in limit of $E$ at infinity is $P_\a$ (i.e., there exists a sequence of numbers $r_n\to \infty$, and the sequence of sets $r_n^{-1}(E)$ converge to $P_\a$ under the Hausdorff distance as $n\to\infty$), then $E$ is a $\P_\a$ set.
\end{thmA}

Besides the global regularity, the property (1.13) helps also to control the the relative distances $d_{x,r}$ between a minimal set and minimal cones in the balls $B(x,r)$ and the local speed of decay of the density function $\theta_x(r)$, because this property gives a lower bound of $\theta_x(r)$. When we prove (1.13) for a minimal cone $C$, we can get nicer local regularity results, that is, if a minimal set is very near $C$ in a ball, then it should be equivalent to $C$ in a smaller ball through a bi-H\"older homeomorphism ($C^1$ diffeomorphism in good cases). So here Theorem \ref{main} has another useful corollary:

\begin{thmB}Let $\theta_1$ be as in Theorem \ref{main}. Then there exists a $\e>0$ such that for any $\a=(\a_1,\a_2)$ with $\a_2\ge \a_1\ge\theta_1$, if $E$ is a 2-dimensional reduced Almgren minimal set in $U\subset\R^4$, $B(x,100r)\subset U$, and there is a reduced minimal cone $P_\a+x$ of type $\P_\a$ centered at $x$ such that $d_{x,100r}(E,P_\a)\le\e$, then there exists a minimal cone $P_{\a'}$ of type $2\P$ such that there is a $C^1$ diffeomorphism $\Phi: B(x,2 r)\to \Phi(B(x,2r))$, such that $|\Phi(y)-y|\le 10^{-2}r$ for $y\in B(x,2r)$, and $E\cap B(x,r)=\Phi( P_{\a'})\cap B(x,r)$.
\end{thmB}

The proof of Theorem \ref{main} will keep us busy until the end of Section 6, but let us already try to explain how it goes.

First notice that Theorem \ref{main} is invariant under translation with respect to $x$, and homogenous with respect to $r$, so we can only restrict to the case where $x=0$ and $r=1$.

Section 2 is devoted to giving some regularity properties for a minimal set $E$ that is close to $P_\a$, but does not contain any point of type $2\P$. In particular, we use a stopping time argument to find a critical region, outside of which everything goes fine, and inside of which things begin to go bad. Here ``bad" means that the set begins to get far away from $P_\a$. The main idea is to control the measure of $E$ in the good region by finer estimates, since there we have good regularity properties; and for the bad region we only control its measure roughly by projections. Part of the argument will be similar to the proof of minimality of $P_\a$.

Section 3 is quite short, where we sum up a little what happens, and give a competitor for $E$, using minimal graphs.

Section 4 and 5 are devoted to giving some estimates for minimal graphs, using some basic estimates for elliptic systems. This leads to some useful control on the measure of the competitor defined in Section 3. 

In Section 6 we conclude, using harmonic extensions and projection properties of the competitor.

We discuss the global regularity and local $C^1$ regularity of minimal sets that are near a $P_\a$ cone in Section 7.

In this article, some of the results and arguments cited in \cite{DJT} exist also in some other (earlier) references, e.g. \cite{Ta}. But for simplify the article, the author will cite \cite{DJT} systematically throughout this article.
%

\noindent\textbf{Some useful notation}

In all that follows, minimal set means Almgren minimal set;

$[a,b]$ is the line segment with end points $a$ and $b$;

$[a,b)$ is the half line with initial point $a$ and passing through $b$;

$B(x,r)$ is the open ball with radius $r$ and centered on $x$;

$\overline B(x,r)$ is the closed ball with radius $r$ and center $x$;

$\overrightarrow{ab}$ is the vector $b-a$;

$H^d$ is the Hausdorff measure of dimension $d$ ;

$d_H(E,F)=\max\{\sup\{d(y,F):y\in E,\sup\{d(y,E):y\in F\}\}$ is the Hausdorff distance between two sets $E$ and $F$.

$d_{x,r}$ : the relative distance with respect to the ball $B(x,r)$, is defined by
$$ d_{x,r}(E,F)=\frac 1r\max\{\sup\{d(y,F):y\in E\cap B(x,r)\},\sup\{d(y,E):y\in F\cap B(x,r)\}\}.$$

\section{A stopping time argument, and regularity and projection properties for minimal sets near $P_\a$}

In this section we use a stopping time argument to control some large scale behavior for minimal sets that near $P_\a$. Let us first introduce some notation.

For each $\a=(\a_1,\a_2)\in [0,\frac\pi 2]^2$ and $i=1,2$, denote by $P_\a=P_\a^1\cup P_\a^2$, where $P_\a^1$ and $P_\a^2$ are two planes in $\R^4$ with characteristic angles $(\a_1,\a_2)$ (this is equivalent to say that there exists an orthonormal basis $\{e_i\}_{1\le i\le 4}$ of $\R^4$ such that $P_\a^1$ is generated by $e_1$ and $e_2$, and $P_\a^2$ is generated by $\cos\a_1e_1+\sin\a_1 e_3$ and $\cos\a_2e_2+\sin\a_2 e_4$). Set
\be C_\a^i(x,r)=(p^i_\a)^{-1}(B(0,r)\cap P^i_\a)+x,\ee
where $p_\a^i$ is the orthogonal projection on $P_\a^i$, and
\be D_\a(x,r)=C_\a^1(x,r)\cap C_\a^2(x,r).\ee
So $C_\a^i(x,r)$ is a cylinder and $D_\a(x,r)$ is the intersection of two cylinders. It is not hard to see that $D_\a(x,r)\supset B(x,r)$ and $D_\a(0,1)\cap P_\a=B(0,1)\cap P_\a$. 

We say that two sets $E,F$ are $\e r$ near each other in an open set $U$ if 
\be d_{r,U}(E,F)<\e,\ee
where 
\be d_{r,U}(E,F)=\frac 1r\max\{\sup\{d(y,F):y\in E\cap U\},\sup\{d(y,E):y\in F\cap U\}\}.\ee

We set also
\be \begin{split}&d^\a_{x,r}(E,F)=d_{r,D_\a(x,r)}(E,F)\\
     &=\frac 1r\max\{\sup\{d(y,F):y\in E\cap D_\a(x,r)\},\sup\{d(y,E):y\in F\cap D_\a(x,r)\}\}.\end{split}\ee

\begin{rem}We should be clear about the fact that 
\be d_{r,U}(E,F)\ne\frac 1rd_H(E\cap U,F\cap  U).\ee
To see this, we can take $U=D_\a(x,r)$, and set $E_n=\partial D_\a(x,r+\frac 1n)$ and $F_n=\partial D_\a(x,r-\frac 1n)$. Then we have 
\be d^\a_{x,r}(E_n,F_n)\to 0\ee and \be \frac 1rd_H(E_n\cap D_\a(x,r),F_n\cap D_\a(x,r))=\frac 1rd_H(E_n\cap D_\a(x,r),\emptyset)=\infty.\ee
So $d_{r,U}$ measures rather how the part of one set in the open set $U$ could be approximated by the other set, and vice versa. However we always have
\be d^\a_{x,r}(E,F)\le\frac 1rd_H(E\cap D_\a(x,r),F\cap D_\a(x,r)).\ee
\end{rem}

Now we give the proposition below, obtained by a stopping time argument.

\begin{pro}\label{rk}There exists $\e_0>0$, such that for any $\e<\e_0$, and $\a>\frac\pi3$, if $E$ is a closed reduced set which is minimal in $D_\a(0,1)$, $d^\a_{0,1}(E,P_\a)<\frac{\e}{10}$, and $E$ contains no $2\P$ point in $B(0,\frac{1}{100})$, then there exists $r_E\in]0,\frac12[$ and $o_E\in B(0,12\e)$ such that $E$ is $2\e r_E$ near $P_\a+o_E$ in $D_\a(o_E, 2r_E(1-12\e))$, but not $\e r_E$ near $P_\a+q$ in $D_\a(o_E,r_E)$ for any $q\in\R^4$.
\end{pro}

\begin{rem}We will also use the construction for information about intermediate scales in the proof. \end{rem}

\noindent Proof of Proposition 2.11.

We fix any $\e$ and $\a=(\a_1,\a_2)>\frac\pi 3$, and set $s_i=2^{-i}$ for $i\ge 0$. Set $D(x,r)=D_\a(x,r),d_{x,r}=d_{x,r}^\a$ for short. 

We proceed in the following way.

Step 1: Denote by $q_0=q_1=O$, then in $D(q_0,s_0)$, $E$ is $\e s_0$ near $P_\a+q_1$ by hypothesis.

Step 2: If in $D(q_1,s_1)$, the set $E$ is not $\e s_1$ near $P_\a+q$ for any $q$, we stop; if not, there exists a $q_2$ such that $E$ is $\e s_1$ near $P_\a+q_2$ in $D(q_1,s_1)$. Here we also ask $\e$ to be small enough (say, $\e<\frac{1}{100}$) so that $q_2\in D(q_1,\frac12s_1)$, thanks to the conclusion of step 1. Then in $D(q_1,s_1)$, we have simultaneously : 
\be d_{q_1,s_1}(E,P_\a+q_1)\le s_1^{-1} d_{q_0,s_0}(E,P_\a+q_1)\le 2\e\ ;\ d_{q_1,s_1}(E, P_\a+q_2)\le\e.\ee

Let us verify that (2.13) implies that $d_{q_1,\frac12 s_1}(P_\a+q_1,P_\a+q_2)\le 12\e$ when $\e$ is small, say, $\e<\frac {1}{100}$. 
In fact, for each $z\in D(q_1,\frac12 s_1)\cap (P_\a+q_1)$, we have $d(z,E)\le d_{q_0,s_0}(E,P_\a+q_1)\le\e$, hence there exists $y\in E$ such that $d(z,y)\le\e$. But since $z\in D(q_1,\frac12 s_1)$, we have $y\in D(q_1,\frac12 s_1+\e)\subset D(q_1,s_1)$, and hence $d(y,P_\a+q_2)\le s_1^{-1} d_{q_1,s_1}(E,P_\a+q_2)\le 2\e$, therefore $d(z,P_\a+q_2)\le d(z,y)+d(y,P_\a+q_2)\le 3\e$.

On the other hand, suppose $z\in D(q_1,\frac12 s_1)\cap (P_\a+q_2)$, we have $d(z,E)\le s_1^{-1}d_{q_1,s_1}(P_\a+q_2,E)\le 2\e$, hence there exists $y\in E$ such that $d(z,y)\le 2\e$. But since  $z\in D(q_1,\frac12 s_1)$, we have $y\in D(q_1,\frac12 s_1+2\e)\subset D(q_0,s_0)$, and hence $d(y,P_\a+q_1)\le d_{q_0,s_0}(E,P_\a+q_1)\le \e$, which implies $d(z,P_\a+q_1)\le d(z,y)+d(y,P_\a+q_1)\le 3\e$.

As a result 
\be d_{q_1,\frac12 s_1}(P_\a+q_1,P_\a+q_2)\le (\frac12 s_1)^{-1}\times 3\e=12\e,\ee
hence $d_{q_1,\frac12 s_1}(q_1,q_2)\le 24\e$, and therefore $d(q_1,q_2)\le 6\e=12\e s_1$.

Now we define our iteration process (notice that it depends on $\e$, so we also call it a $\e$-process).

Suppose that all $\{q_i\}_{i\le n}$ have already been defined, with
\be d(q_i,q_{i+1})\le 12s_i\e=12\times 2^{-i}\e\ee
for $0\le i\le n-1$, and hence
\be d(q_i,q_j)\le 24\e s_{\min(i,j)}=2^{-\min(i,j)}\times 24\e\ee 
for $0\le i,j\le n$. Moreover, for all $i\le n-1$, $E$ is $\e s_i$ near $P_\a+q_{i+1}$ in $D(q_i, s_i)$. We say that the process does not stop at step $n$. In this case

Step n+1 : We look at the situation in $D(q_n, s_n)$.

If $E$ is not $\e$ near any $P_\a+q$ in this ball of radius $s_n$, we stop, since we have found the $o_k=q_n, r_k=s_n$ as desired. In fact, since $d(q_{n-1}, q_n)\le 12\e s_{n-1}$, we have
$D(q_n, 2s_n(1-12\e))=D(q_n, s_{n-1}(1-12\e))\subset D(q_{n-1}, s_{n-1})$, and hence
\be\begin{split}
d_{q_n,2s_n(1-12\e)}(P_\a+q_n, E)&\le (1-12\e)^{-1}d_{q_{n-1},s_{n-1}}(P_\a+q_n, E)\\
&\le \frac{\e}{1-12\e}.
\end{split}\ee
Moreover 
\be d(o_k, O)=d(q_n,q_1)\le 2^{-\min(1,n)}\times 24\e=12\e.\ee

Otherwise, we can find a $q_{n+1}\in \R^4$ such that $E$ is still $\e s_n$ near $P_\a+q_{n+1}$ in $D(q_n, s_n)$. Then since $\e$ is small, $q_{n+1}\in D(q_n,\frac12 s_n)$. Moreover we have as before $d(q_{n+1},q_n)\le 12\e s_n$, and for $i\le n-1$,
\be d(q_i,q_{n+1})\le \sum_{j=i}^n d(q_j,q_{j+1})\le\sum_{j=i}^n 12\times 2^{-j}\e\le 2^{-j}\times 24\e=2^{-\min(i,n+1)}\times 24\e.\ee
Thus we have obtained our $q_{n+1}$.  

Now all we have to do is to prove that for every $\e$ small enough, this process has to stop at a finite step. For this purpose we need the following proposition.

\begin{pro} There exists $\theta_1'\in[\theta_0,\frac\pi 2)$, and for any $l\in]0,\frac12]$, there exists $\e_l\in]0,\frac 12[$, such that for any $\a>\theta_1'$, $\e\le \e_l$, and $E$ as in Proposition 2.11, if the $\e-$process does not stop before the step $n$, then

(1) The part $E\cap (D_\a(0,\frac {39}{40})\backslash D_\a(q_n,\frac {1}{10}s_n))$ is composed of two disjoint pieces $G^i,i=1,2$, such that:
\be
G^i\mbox{ is the graph of a }C^1\mbox{ map }\ g^i:C^i_\a(0,\frac{39}{40})\backslash C^i_\a(q_n,\frac{1}{10}s_n)\cap P^i_\a\to {P^i_\a}^\perp\ee
with
\be||\nabla g^i||_\infty<l\le\frac12;\ee

(2) For every $\frac {1}{10}s_n\le t\le s_n$
\be E\cap (D_\a(0,1)\backslash D_\a(q_n,t))=G_t^1\cup G_t^2\ee
where $G_t^1, G_t^2$ do not meet each other. Moreover
\be P_\a^i\cap (D_\a(0,1)\backslash C_\a^i(q_n,t))\subset p_\a^i(G_t^i)\ee
where $p_\a^i$ is the orthogonal projection on $P_\a^i,i=1,2$;

%
\end{pro}

\begin{rem}If we take the optimal $\e_l$ for each $l$ such that Proposition 2.20 holds, then obviously for any $l\le l'$, $\e_{l}\le\e_{l'}$.
\end{rem}

We will not prove this proposition, see \cite{2p} Proposition 6.1 (1) (2) for the proof. But we'll use it to finish our Proposition 2.11.

\begin{rem} In fact we need all the properties stated in \cite{2p} Proposition 6.1 for our set $E$. For (1) and (2) in \cite{2p} Proposition 6.1, the arguments there can be applied directly here to our set $E$ with no change. But for (3) and (4), the proof in \cite{2p} Proposition 6.1 uses some special property of $E_k$, which are not necessarily true for our set $E$ here. Hence we will treat  the property of surjective projections ( (4) of \cite{2p} Proposition 6.1) later in a different way.
\end{rem}

So let $\e_0$ be the $\e_\frac 12$ in Proposition 2.20. Suppose that the $\e-$process does not stop at any finite step, and we'll try to get a contradiction. By (1) of Proposition 2.20, for any $n$, $E\cap (D_\a(0,1)\backslash D_\a(q_n,\frac{1}{10}s_n))$ is composed of two disjoint graphs $G^i$ on $[C^i_\a(0,1)\backslash C^i_\a(q_n,\frac{1}{10}s_n)]\cap P_\a^i, i=1,2$. Denote by $\Delta_n=D_\a(q_n,s_n)$.

Notice that by (2.19), with $\e<\frac{1}{100}$, the sets $\Delta_n=D_\a(q_n,s_n)$ are in fact a sequence of non degenerate compact balls, with
\be  \Delta_n\subset \Delta_{n-1}, n\in\N, \lim_{n\to\infty}\mbox{diam}(\Delta_n)\to 0,\ee
Hence there exists a point $p\in B(0,\frac 12)$, such that $\{p\}=\cap_n \Delta_n$. Then $p$ is also the limit of $q_n$, hence it lies in $B(0,\frac{1}{100})$. By (1) of Proposition 2.20, for any $r\in(0,\frac 12)$, $E\cap D(p,\frac 12)\bs D(p,r)$ is composed of the union of two disjoint graphs on $P_\a^i\cap C_\a^i(p,\frac 12)\bs C_\a^i(p,r)$. As a result , $E\cap D(p,\frac12)\bs\{p\}$ is composed of two $C^1$ graphs on $P_\a^i\cap C_\a^i(p,\frac 12)\bs \{p\}$. Denote by $G^i$ these two graphs. By (2.22), they are both $\frac12$-Lipschitz. Now $E$ is closed hence $p\in E$. Then for each $i=1,2$, $G^i\cup\{p\}$ is a $\frac 12$-Lipschitz graph on $P_\a^i\cap C_\a^i(p,\frac 12)$, and hence $E\cap D_\a(p,\frac 12)$ is composed of the disjoint union of these two $\frac 12$-Lipschitz graphs. Now we define $\varphi:E\cap D_\a(p,\frac 12)\to P_\a+p$, where the restriction of $\varphi$ to each $G^i\cup\{p\}$ is just the orthogonal projection to $P_\a^i+p$. Then it is easy to check that $\varphi$ is a Lipschitz homeomorphism. That is, $E$ is bi-Lipschitz homeomorphic to $P_\a$ in $D_\a(p,\frac 12)$. 

We want to prove that $p$ is a point of type $2\P$. Take any blow-up limit $C$ of $E$ at the point $p$. Then $C$ is a minimal cone. By the bi-H\"older regularity for 2-dimensional minimal sets, near the point $p$, $E$ is locally bi-H\"older equivalent to $C$. But $E$ is also bi-Lipschitz equivalent to $p_\a$ near $p$, hence the two minimal cones $P_\a$ and $C$ are topologically the same. As a consequence, $P_\a\cap\partial B(0,1)$ and $C\cap \partial B(0,1)$ are topologically the same, therefore, $C\cap \partial B(0,1)$ is the union of two topological circles. But by the description of 2-dimensional minimal cones (cf.\cite{DJT}, Proposition 14.1), the intersection of any minimal cone with the unit sphere is a finite union of great circles and arcs of great circles that meet at their extremities by group of three with $120^\circ$ angles. Here in our case, we can deduce that $C\cap\partial B(0,1)$ is the union of two circles. Hence $C$ is a minimal cone of type $2\P$.

Hence the point $p$ is a point of type $2\P$. This contradicts the fact that $E\cap B(0,\frac{1}{100})$ contains no point of type $2\P$, because $p\in B(0,\frac{1}{100})$.

Thus we complete the proof of Proposition 2.11.\qed

Next we still have to prove some property of surjective projection, as remarked in Remark 2.26. 

\begin{pro} Take $\e\le\e_0$, and take $\a$ and $E$ as in Proposition 2.20. Then for any $n\ge 1$, if the $\e-$process does not stop before the step $n$, then the orthogonal projections $p_\a^i: E\cap \overline D_\a(q_n,t)\to P_\a^i\cap \overline C_\a^i(q_n,t),i=1,2$ are surjective, for all $\frac 19s_n\le t\le s_n$.
\end{pro}

\nd Fix a such $n$. Set $s_i=2^{-i}$ for $i\ge 0$. Set $D(x,r)=D_\a(x,r),C^i(x,r)=C^i_\a(x,r), d_{x,r}=d_{x,r}^\a$ for short. By (1) of Proposition 2.20,  the part $E\cap (D_\a(0,\frac {39}{40})\backslash D_\a(q_n,\frac {1}{10}s_n))$ is composed of two disjoint pieces $G^i,i=1,2$, such that:
\be
G^i\mbox{ is the graph of a }C^1\mbox{ map }\ g^i:C^i_\a(0,\frac{39}{40})\backslash C^i_\a(q_n,\frac{1}{10}s_n)\cap P^i_\a\to {P^i_\a}^\perp\ee
with
\be||\nabla g^i||_\infty<\frac12.\ee

Thus $G^i\cap\partial C^i(0,\frac {39}{40})$ is a nice $C^1$ curve, which is the graph of $g^i$ on $P_\a^i\cap \partial C^i(0,\frac {39}{40})$, and $g^i$ is $\frac 12$-Lipschitz. Denote by $\gamma^i=g^i|_{P_\a^i\cap \partial C^i(0,\frac {39}{40})}$. Then $||\gamma^i||_\infty\le \frac{\e}{10}$ by hypothesis.

Now we define a set $Q$ as follows. First, $Q\subset\overline B(0,1)$, and $Q\bs D(0,\frac {39}{40})=E\bs D(0,\frac {39}{40})$. Inside $D(0,\frac 34)$, $Q\cap \overline D(0,\frac 34)=P_\a\cap \overline D(0,\frac 34)$, the union of two planes. For the part on the annulus $D(0,\frac {39}{40})\bs \overline D(0,\frac 34)$, we just use two graphs of affine functions to join $P_\a^i\cap \partial D(0,\frac 34)$ and $\gamma^i$. That is, we define $h^i: P_\a^i\cap D(0,\frac {39}{40})\bs \overline D(0,\frac 34)\to {P_\a^i}^\perp$, for any $x\in P_\a^i\cap D(0,\frac {39}{40})\bs \overline D(0,\frac 34) (\frac 34,\frac {39}{40})$, $h^i(x)=\frac{|x|-\frac 34}{\frac{39}{40}-\frac 34} \gamma^i(\frac{39x}{40|x|})$.

Thus for any $x\in D(0,\frac {39}{40})\bs \overline D(0,\frac 34)$, $|\frac{\partial}{\partial r}h^i(x)|=\frac{1}{\frac{39}{40}-\frac 34}|\gamma^i(\frac{39x}{40|x|})|\le \frac{40}{9}\frac{\e}{100}\le \frac{\e}{20}\le \frac{1}{2000}$, and $|\frac{\partial}{\partial\theta}(x)|\le Lip(\gamma^i)\le\frac 12$, hence the tangent direction derivative is less than 
\be \frac{1}{|x|}|\frac{\partial}{\partial\theta}(x)|\le\frac12/\frac 34=\frac 23.\ee
Hence we have 
\be \mbox{Lip }h^i\le \max\{\frac{1}{2000},\frac 23\}=\frac 23.\ee
Thus the map $H^i:P_\a^i\cap D(0,\frac {39}{40})\bs \overline D(0,\frac 34)\to\R^4: x\mapsto (x,h^i(x))$ is $(1+(\frac23)^2)^\frac 12=\frac{\sqrt {13}}{3}$-Lipschitz. So if we denote by $\Sigma^i$ the graph of $h^i$, then
\be \begin{split}
H^2(\Sigma^i)&=H^2(H^i(P_\a^i\cap D(0,\frac {39}{40})\bs \overline D(0,\frac 34))\le (\frac{\sqrt {13}}{3})^2)H^2(P_\a^i\cap D(0,\frac {39}{40})\bs \overline D(0,\frac 34))\\
&=\frac{897}{1600}\pi\le\frac{9\pi}{16}, i=1,2.\end{split}\ee

\centerline{\includegraphics[width=0.6\textwidth]{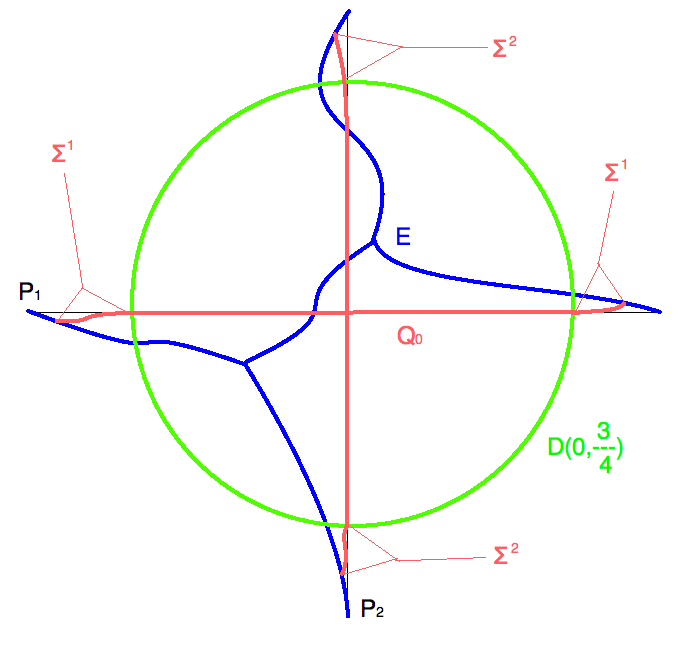}}
\nopagebreak[4]
\centerline{2-1}

Let $Q=[E\bs D(0,\frac {39}{40})]\cup\Sigma^1\cup\Sigma^2\cup[P_\a\cap D(0,\frac 34)]$, and $Q_0=Q\cap D(0,\frac {39}{40})$. (See Figure 2-1.) Set $Q^i=\Sigma^i\cup[P_\a^i\cap D(0,\frac 34)]$, then $Q_0$ is the almost disjoint union $Q^1\cup Q^2$. For each $i=1,2$, 
\be H^2(Q^i)=H^2(\Sigma^i)+H^2(P_\a^i\cap D(0,\frac 34))\le\frac{9\pi}{16}+\frac{9\pi}{16}=\frac{9\pi}{8}.\ee

Notice that the set $Q_0$ is a $C^1$ version of $P_\a\cap D(0,\frac34)$, and $Q^i,i=1,2$ are its two flat parts as $P_\a^i$. 

Now suppose that for some $t\in[\frac19s_n,s_n)$, for example the projection $p_\a^1: E\cap D(q_n,t)\to P_\a^1\cap C^1(q_n,t)$ is not surjective. Then we are going to prove that we can deform $E$ to $[Q\bs Q_0]\cup Q^2$, and deduce a contradiction.

So take a point $p\in P_\a^1\cap \overline C^1(q_n,t)$ which does not admit a pre-image in $E\cap \overline D(q_n,t)$. Since the set $E_t:=E\cap \overline D(q_n,t)$ is compact, its projection $p_\a^1(E_t)$ is also compact, which means that we can pick $p\in P_\a^1\cap C^1(q_n,t)\bs p_\a^1(E_t)$ and $r\in(0,\frac{t}{10})$ such that $B(p,r)\cap P_\a^1\subset P_\a^1\cap C^1(q_n,t)\bs p_\a^1(E_t)$, and moreover $0\not\in B(p,3r)$.

Now the set $E_t\subset\overline D(q_n,t)\bs {p_\a^1}^{-1}(B(p,r)\cap P_\a^1)$. Take an orthogonal union of two planes $P_0=P_0^1\cup_\perp P_0^2$ in $\R^4$, denote by $p_0^i$ the orthogonal projection on $P_0^i,k=1,2$, take a point $p_0\in P_0^1$ such that $d(p_0,o)=\frac12$.

Then we can easily find a Bi-Lipschitz mapping $\varphi:\overline D(q_n,t)\bs {p_\a^1}^{-1}(B(p,r)\cap P_\a^1)\to \overline D(0,1)\bs {p_0^1}^{-1}(B(p_0,\frac14)\cap P_0^1)$, such that $\varphi(E_t\cap D(q_n,t)\bs D(q_n,\frac{1}{10} s_n))=P_0\cap D(0,1)\bs D(0,\frac34)$ (because in the annulus $D(q_n,t)\bs D(q_n,\frac{1}{10} s_n)$, the set $E$ is still a $C^1$ graph of $P_\a$).

For any point $x\in D(0,1)$, write $x=(x_1,x_2)$, where $x_i=p_0^i(x)\in B^i(0,1),i=1,2$ ($B^i(0,1)$ is the unit ball of the plane $P_0^i$). We define $\psi:D(0,1)\bs {p_0^1}^{-1}(B((p_0,\frac14)\cap P_0^1)\to D(0,1)\cap P_0\bs {p_0^1}^{-1}(B((p_0,\frac14)\cap P_0^1)$ as follows:
\be \psi(x)=\left\{\begin{array}{cc}p_0^1(x), & x_2<\frac 34;\\ (x_1,4x_2-3), & x_2\ge\frac34.\end{array}\right.\ee

Then $\psi$ is a Lipschitz map, which maps $[C^1(0,1)\cap C^2(0,\frac34)]\cup [P_0\cap D(0,1)]$ to $P_0\cap D(0,1)$, and $\psi|_{P_0\cap\partial D(0,1)}=Id$. In particular, $\psi(\varphi(E_t))\subset P_0\cap D(0,1)\bs {p_0^1}^{-1}(B(p_0,\frac14)\cap P_0^1)$.

Thus the map $f_1=\varphi^{-1}\circ\psi\circ\varphi$ maps $E_t$ to $P_\a\cap D(q_n,t)\bs D(q_n,\frac{1}{10} s_n)$, and $f_1|_{E\cap\partial D(q_n,t)}=id$.

We can extend $f_1$ to a Lipschitz map from $D(0,\frac{39}{40})\to D(0,\frac{39}{40})$, such that $f_1|_{E\cap D(0,\frac{39}{40})\bs D(q_n,t)}=id$ and $f_1|_{D(0,\frac{39}{40})\bs D(0,\frac 12)}=id$.

Then $f_1$ is a deformation of $E$ in $D(0,\frac{39}{40})$, which sends $E\cap D(0,\frac{39}{40})$ to $Q_0\bs [B(p,r)\cap P_\a^1]$, this is the union of $Q^2$ and $Q^1$ minus a hole $B(p,r)\cap P_\a^1$. So we can keep on the deformation, and take the map $f_2$ which deforms $Q^1\bs [B(p,r)\cap P_\a^1]$ to a set $E^1=\{0\}\cup\partial Q^1\cup C$ of measure zero, where $C$ is a segment that connects the origin and $\partial Q^1$ and keeps $Q^2$ fixed. Then the map $f=f_2\circ f_1$ sends $E\bs D(0,\frac{39}{40})$ to $Q^2\cup E^1$, hence the measure
\be H^2(E\cap D(0,\frac{39}{40}))=H^2(Q^2)\le\frac{9\pi}{8}.\ee

The map $f$ is Lipschitz, and its restriction to $Q_0\cap\partial D(0,\frac{39}{40})$ is the identity. We extend $f$ to a Lipschitz map on $D(0,1)$, still denoted by $f$, such that $f=id$ near the boundary of $D(0,1)$. Thus by the minimality of $E$, and since $f$ does not move $E\bs D(0,\frac{39}{40})$, we have
\be H^2(E\cap D(0,\frac{39}{40}))\le H^2(f(E\cap D(0,\frac{39}{40}))\le\frac{9\pi}{8}.\ee

However since $n>1$, we have $s_n<\frac 12$. By (1) of Proposition 2.20, we have
\be \begin{split}H^2(E\cap D(0,\frac{39}{40}))&\ge H^2(G^1)+H^2(G^2)\ge H^2(p_\a^1(G^1))+H^2(p_\a^2(G^2))\\
&=\sum_{i=1,2}H^2(P_\a^i\cap C^i(0,\frac{39}{40})\bs C^i(q_n,\frac{1}{10}s_n))\\
&\ge \sum_{i=1,2}H^2(P_\a^i\cap C^i(0,\frac{39}{40})\bs C^i(q_n,\frac{1}{20}))\\
&=2\times \pi((\frac{39}{40})^2-(\frac{1}{20})^2)=\frac{1517}{800}\pi>\frac{9\pi}{8},\end{split}\ee
which leads to a contradiction.

This completes the proof of Proposition 2.28.\qed

\section{A competitor, and estimates for minimal graphs}

Let $\theta_1'$, $\a$ be as in Proposition 2.20, let $ \e=\e_0$, $\mu$ be chosen later, and let $E$ be as in Proposition 2.11, that is, $d_{0,1}^\a<\frac{\e}{10}$, and $E$ contains no $2\P$ type point in $B(0,\frac{1}{100})$. We want to construct a competitor for $E$, and show that if $d_{0,1}^\a$ is sufficiently small, this competitor admits necessarily less measure than $E$, and thus leads to a contradiction.

Let us point out that the condition $d_{0,1}^\a<\frac{\e}{10}$ is a general qualitative one, which guarantees that $E$ satisfies the  regularity properties in Proposition 2.20 and 2.28. To make the necessary finer estimates for measures of $E$ and its competitor, we still have to get the "$\lambda$-near" condition as in Theorem \ref{main}. 

\medskip

So by Proposition 2.11, there is a $r_E\in]0,\frac12[$, $o_E\in B(0,\frac12\e_0)$ such that the conclusion in Proposition 2.11 holds for $E$. Denote by $\gamma^i:\partial B(0,\frac 12)\cap P_\a^i\to {P_\a^i}^\perp$ the $C^1$ curve $g^i|_{\partial B(0,\frac 12)\cap P_\a^i}$. Suppose that $||\gamma^i|_{\partial B(0,\frac 12)\cap P_\a^i}||_{C^1}\le\mu$.

The idea of the construction of the competitor is not complicated. We take, for each $i$, a minimal graph $\Sigma^i$ which is the graph of a function $f^i:B(0,\frac 12)\cap P_\a^i\to {P_\a^i}^\perp$ such that $f^i|_{\partial B(0,\frac 12)\cap P_\a^i}=\gamma^i$. Take $\Sigma=\Sigma^1\cup\Sigma^2$. Then hopefully when $\mu$ is small enough, these two graphs are very flat at the center, so that $\Sigma$ is very similar to $P_\a$. Thus we can deform $E\cap D_\a(0,\frac 12)$ to a subset of $\Sigma$ in a Lipschitz manner, while keeping $E\cap\partial D_\a(0,\frac 12)$ unchanged. Hence $\Sigma$ contains a competitor of $E$ in $D_\a(0,\frac 12)$. By the minimality of $E$, the measure of $\Sigma$ has to be larger than that of $E\cap D_\a(0,1)$. But we are going to show that when $\mu$ is small enough, this is not true.

Before we go down to the following two sections, which will be devoted to giving some estimates for minimal graphs, let us already explain what happens. 

We want to compare the measures of $E\cap D_\a(0,\frac 12)$ and $\Sigma$. Outside $D(o_E,\frac{1}{10}r_E)$, by Proposition 2.20, $E$ is also composed of two $C^1$ graphs $G^i$ on the two annuli $P_\a^i\cap B(0,\frac 12)\bs C^i((o_E,\frac{1}{10}r_E)$. So in this part, our goal is to compare the surface measure of $\Sigma^i$ and $G^i$, that is, the graph of $f^i$ and $g^i$. Notice that $f^i$ and $g^i$ coincide on $\partial P_\a^i\cap \partial B(0,\frac 12)$, and on $P_\a^i\cap\partial B(o_E,\frac{1}{10}r_E)$, $g^i$ is supposed to be $\e$-far from any plane, while $f^i$ is almost a plane (this is the main result of Section 4). Then Section 5 is devoted to estimating the difference between these two graphs.

So this will help estimate the difference between measures of $E$ and $\Sigma$ on the annulus region $D_\a(0,\frac12)\bs D(o_E,\frac{1}{10}r_E)$. For the part of $E\cap D(o_E,\frac{1}{10}r_E)$, we estimate its measure by using projections.

\section{Existence and estimates for derivatives for minimal graphs}

Denote by $B=B(0,1)\cap \R^2$ the unit disc in $\R^2$. Let $\gamma$ be a $C^1$ function from $\partial B$ to $\R^2$. Now by Theorems 4.1 and 4.2 of \cite{LaOs}, there exists a function $f:\overline B\to \R^2$, whose graph $\Sigma_f=\{(x,f(x)):x\in\overline B\}\subset\R^4$ is a minimal surface, $f|_{\partial B}=\gamma$, and $f\in C^0(\overline B)\cap C^\infty(B).$ In particular, by (c) of Theorem 4.1 of \cite{LaOs} and the maximum principle for harmonic maps, we have 
\be||f||_\infty\le||\gamma||_{L^\infty(\partial B)}.\ee

Now suppose that $\mu=\max\{||\gamma||_{L^\infty(\partial B)},||D\gamma||_{L^\infty(\partial B)}\}$ is small, then by (4.1), $||f||_\infty\le\mu$ is small. We want to prove that $|\nabla f|,|\nabla^2 f|,|\nabla^3f|$ are also small in a neighborhood of 0, and are controlled by $\mu$. More precisely, we state the following proposition.

\begin{pro}There exists $\mu_0>0$, such that for any $\mu<\mu_0$, there exists a constant $C(\mu)$, with $\lim_{\mu\to 0}C(\mu)=0$, such that if $f$ is a minimal graph on $B(0,1)$, with 
\be\max\{||f|_{\partial B(0,1)}||_\infty, ||Df|_{\partial B(0,1)}||_\infty\}\le\mu,\ee then
\be \max_{0\le i\le 3}||\nabla^i f||_{L^\infty(B(0,\frac 34))}\le C(\mu).\ee
\end{pro}

\nd

First let us apply a regularity theorem on varifolds to get the initial estimate for $\nabla f$, and then we can go into the machine of estimates for elliptic systems. Before stating the theorem, we give some useful notations below.

 $G(n,d)$ denotes the Grassmann manifold $G(\R^n,d)$; 
 
for every $T\in G(n,d)$, we denote by $\pi_T$ the orthogonal projection on the $d$-plane represented by $T$;
 
 for every measure $\nu$ on $\R^n$, $\theta^d(\nu,x)=\lim_{r\to 0}\frac{\nu B(a,r)}{\a(d)r^d}$ (if the limit exists) is the density of $\nu$ on $x$, where $\a(d)$ denotes the volume of the $d$-dimensional unit ball; 
 
 $\mathbb V_d(\R^n)$ denotes the set of all $d-$varifold in $\R^n$, i.e. all Radon measures on $G_d(\R^n)=\R^n\times G(n,d)$; 
 
 for each $V\in \mathbb V_d(\R^n),$ $||V||$ is the Radon measure on $\R^n$ such that for each $A\subset \R^n$, $||V||(A)=V(G_d(\R^n)\cap \{(x,S):x\in A\})$; 
 
 $\delta(V)$ denotes the first variation of $V$, that is, the linear map from $\mathfrak X(R^n)$ to $\R$, defined by
\be\delta V(g)=\int Dg(x)\cdot \pi_S dV(x,S)\ee for $g\in\mathfrak X(\R^n).$
Here $\mathfrak X (\R^n)$ is the vector space of all $C^\infty$ maps from $\R^n$ to $\R^n$ with compact support.

In our case, we are only interested in rectifiable varifolds. In fact, with each $d-$rectifiable set $E$ we associate a $d-$varifold, denoted by $V_E$, in the following sense: for each $B\subset \R^n\times G(n,d)$, we have 
\be V_E(B)=H^d\{x:(x,T_xE)\in B\}.\ee
Recall that $T_xE$ is the $d$-dimensional tangent plane of $E$ at $x$; it exists for almost all $x\in E$, because $E$ is $d-$rectifiable. Then $||V_E||=H^d|_E$. Moreover, the density $\theta^d(||V_E||,x)$ exists for almost all $x\in E$.

\begin{thm}[cf.\cite{All72} Regularity theorem at the beginning of section 8]Suppose $2\le d<p<\infty,$ $q=\frac{p}{p-1}$. Corresponding to every $\e\in]0,1[$ there is $\eta>0$ with the following property:

Suppose $0<R<\infty$, $0<\lambda<\infty$, $V\in\mathbb{V}_d(\R^n)$, $a\in spt||V||$ and

1) $\theta^d(||V||,x)\ge \lambda$ for $||V||$ almost all $x\in B(a,R)$;

2) $||V||B(a,R)\le(1+\eta)\lambda\a(d)R^d ;$

3) $\delta V(g)\le\eta \lambda^{\frac 1p}R^\frac{d}{p-1}(\int|g|^q\lambda||V||)^\frac 1q$ whenever $g\in {\mathfrak X}(\R^n)$ and $spt\ g\subset B(a,R)$.

Then there are $T\in G(n,d)$ and a continuously differentiable function $F:T\to \R^n$, such that $\pi_T\circ F=1_T$,
\be||DF(y)-DF(z)||\le\e(|y-z|/R)^{1-\frac dp}\mbox{ whenever }y,z\in T,\ee and
\be B(a,(1-\e)R)\cap spt||V||=B(a,(1-\e)R)\cap\ image\ F.\ee
\end{thm}

\begin{rem}

1) In the theorem, since $\pi_T\circ F=1_T$,  we can see that $F$ is in fact the graph of a $C^1$ function $f$, defined by $f(t)=\pi_{T^\perp} F(t)$, with $t\in T$, $\pi_{T^\perp}$ the orthogonal projection on the orthogonal space $T^\perp$ of $T$. Moreover $||Df(t)||\le ||DF(t)||$ for all $t\in T$.

2) If $E$ is a minimal surface, then $V_E$ is stationary, i.e. $\delta V_E=0$. Hence the condition 3) is automatically true. In fact if we set $g_t(x)=(1-t)x+tg(x)$, then 
\be\delta V_E(g)=\frac{d}{dt} H^d(g_t(E\cap spt g)),\ee 
which can be deduced from the area formula. Thus if $E$ is a minimal surface, $\delta V_E=0$. 
\end{rem}

Now we want to apply Theorem 4.7 to our set $\Sigma_f$, so we have to check all the conditions in the theorem. We take $\lambda=1,a=(0,f(0)), R=1$, then 1) is true, by the fact that $\Sigma_f$ is a $C^\infty$ manifold; 3) is true by the Remark 4.10 2); for 2), notice first of all that $B(a,R)\cap \Sigma_f\subset\Sigma_f$, so we just have to estimate the surface of $\Sigma_f$. Notice that Lip  $\gamma\le \mu$, hence for the length of the graph of $\gamma$, denoted also by $\gamma$, we have
\be |\gamma|=\int_{\partial B}\sqrt{1+|D\gamma|^2}\le \int_{\partial B}\sqrt{1+\mu^2}=2\pi(1+\mu^2).\ee

Now by the isoperimetric inequality for minimal surface (cf. \cite{Car}), we have
\be 4\pi H^2(\Sigma_f)\le |\gamma|^2=[2\pi(1+\mu^2)]^2,\ee
which means
\be H^2(\Sigma_f\cap B(a,R))\le H^2(\Sigma_f)\le(1+\mu^2)^2\pi.\ee

Hence we can take $\mu$ small enough such that 2) holds for some $\eta$, such that (4.8) and (4.9) are true for some $\e$ small, which give us that
\be ||f||_{C^{1,\sigma}(B(0,\frac 89))}\le C_1(\mu),\ee
with $\lim_{\mu\to 0}C_1(\mu)=0.$

\begin{rem}We might be able to use only the estimates for elliptic system to get this initial estimate, without using the powerful Theorem 4.7.  
\end{rem}

For estimating higher order regularity of $f$, we have to refer to the minimal surface equation system and put everything in the machine of elliptic system. 

First we give some notations. 

Denote by $M_2(\R)$ the set of $2\times 2$ matrices on $\R$. For any $\left(\begin{array}{cc}a\ & b\\c\ &d\end{array}\right)\in M_2(\R)$, denote by $|\left(\begin{array}{cc}a\ & b\\c\ &d\end{array}\right)|=a^2+b^2+c^2+d^2$, and for any $\left(\begin{array}{cc}a'\ & b'\\c'\ &d'\end{array}\right)\in M_2(\R)$, define $<\left(\begin{array}{cc}a\ & b\\c\ &d\end{array}\right),\left(\begin{array}{cc}a'\ & b'\\c'\ &d'\end{array}\right)>=aa'+bb'+cc'+dd'.$ Denote by $\cdot$ the multiplication of matrices. Set, for any $\left(\begin{array}{cc}a\ & b\\c\ &d\end{array}\right)\in M_2(\R)$, $\left(\begin{array}{cc}a\ & b\\c\ &d\end{array}\right)^*=\left(\begin{array}{cc}d\ & -c\\-b\ &a\end{array}\right)\in M_2(\R).$

\smallskip

For any domain $\Omega\subset\R^2$, for any differentiable function $h:\Omega\to\R$, denote by $h_x,h_y$ its two partial derivates. For any $C^2$ function $h=(h^1,h^2):\R^2\to\R^2$, with $h^i:\R^2\to\R$ two $C^1$ functions, denote by $\nabla h$ the matrix valued function $\left(\begin{array}{cc}h^1_x\ & h^2_x\\h^1_y\ &h^2_y\end{array}\right)$. And for any matrix valued function $f=\left(\begin{array}{cc}f^1\ & f^2\\f^3\ &f^4\end{array}\right)$ on $\R^2$, we define $div f=(f^1_x+f^3_y,f^2_x+f^4_y)\in\R^2$.

Then we have
\be H^2(\Sigma_h)=\int_\Omega\sqrt{1+|\nabla h|^2+(\det \nabla h)^2}.\ee

Denote by $S(h)=|\nabla h|^2+(\det \nabla h)^2$ for any $h$. 

$\Sigma_f$ is a minimal submanifold, hence it is stable with respect to any local perturbation. More precisely, for any $C^\infty$ function $\varphi:\overline B\to\R^2$ with $\varphi|_{\partial B}=0_{\R^2}$, we have
\be \frac{d}{dt}|_{t=0}H^2(\Sigma_{f+t\varphi})=0.\ee
(4.17) and (4.18) gives that, for any $C^\infty$ function $\varphi:\overline B\to\R^2$ with $\varphi|_{\partial B}=0_{\R^2}$,
\be \begin{split}0&=\frac{d}{dt}|_{t=0}\int_B\sqrt{1+|\nabla (f+t\varphi)|^2+(\det \nabla (f+t\varphi))^2}\\
&=\int_B\frac{d}{dt}|_{t=0}\sqrt{1+|\nabla (f+t\varphi)|^2+(\det \nabla (f+t\varphi))^2}\\
&=\int_B\frac{\frac{d}{dt}|_{t=0}<\nabla (f+t\varphi),\nabla (f+t\varphi)>+\frac{d}{dt}|_{t=0}(\det \nabla (f+t\varphi))^2}{2\sqrt{1+S(f)}}\\
&=\int_B\frac{<\nabla f,\nabla \varphi>+\det(\nabla f)\frac{d}{dt}|_{t=0}(\det \nabla (f+t\varphi))}{\sqrt{1+S(f)}}.\end{split}\ee

Denote by $\nabla f=\left(\begin{array}{cc}A\ & B\\C\ &D\end{array}\right)$, and $\nabla \varphi=\left(\begin{array}{cc}a\ & b\\c\ &d\end{array}\right)$, then we have
\be \begin{split}\det \nabla (f+t\varphi)&=\det\left(\begin{array}{cc}A+ta\ & B+tb\\C+tc\ &D+td\end{array}\right)\\
&=(A+ta)(D+td)-(B+tb)(C+tc)\\
&=\det\nabla f+t^2\det\nabla\varphi+t(aD-bC-cB+dA)\\
&=\det\nabla f+t^2\det\nabla\varphi+t<(\nabla f)^*,\nabla\varphi>.\end{split}\ee
Therefore
\be \frac{d}{dt}\Big|_{t=0}(\det \nabla (f+t\varphi))=<(\nabla f)^*,\nabla\varphi>.\ee

Combining with (4.19), we get
\be \int_B<\frac{\nabla f+det(\nabla f)(\nabla f)^*}{\sqrt{1+S(f)}},\nabla\varphi>=0\ee
for any $C^\infty$ function $\varphi:\overline B\to\R^2$ with $\varphi|_{\partial B}=0$. Hence we have
\be div(\frac{\nabla f+det(\nabla f)(\nabla f)^*}{\sqrt{1+S(f)}})=(0,0).\ee

This means, $f$ satisfies the elliptic system (4.23). Denote by $f=(u,v)$, with $u,v$ two functions from $\overline B$ to $\R$. Denote by $u_x,u_y,v_x,v_y$ the partial derivatives of $f$ for short, and we write the system (4.23) in the standard non-linear form below
\be \left\{\begin{split}
\frac{\partial}{\partial x}[\frac{(1+v_y^2)u_x-(v_xv_y)u_y}{\sqrt{1+S(f)}}]+\frac{\partial}{\partial y}[\frac{(1+v_x^2)u_y-(v_xv_y)u_x}{\sqrt{1+S(f)}}]&=0,\\
\frac{\partial}{\partial x}[\frac{(1+u_y^2)v_x-(u_xu_y)v_y}{\sqrt{1+S(f)}}]+\frac{\partial}{\partial y}[\frac{(1+u_x^2)v_y-(u_xu_y)v_x}{\sqrt{1+S(f)}}]&=0.\end{split}\right.
\ee

%

Now set, for any $\left(\begin{array}{cc}a&b\\c&d\end{array}\right)\in M_2(\R)$, $T\left(\begin{array}{cc}a&b\\c&d\end{array}\right)=1+a^2+b^2+c^2+d^2+(ad-bc)^2$, and
\be \begin{split}A_x^x\left(\begin{array}{cc}a&b\\c&d\end{array}\right)=\frac{(1+d^2)a-bcd}{\sqrt{T\left(\begin{array}{cc}a&b\\c&d\end{array}\right)}}, 
A_y^x\left(\begin{array}{cc}a&b\\c&d\end{array}\right)=\frac{(1+b^2)c-abd}{\sqrt{T\left(\begin{array}{cc}a&b\\c&d\end{array}\right)}}, \\
A_x^y\left(\begin{array}{cc}a&b\\c&d\end{array}\right)=\frac{(1+c^2)b-acd}{\sqrt{T\left(\begin{array}{cc}a&b\\c&d\end{array}\right)}}, 
A_y^y\left(\begin{array}{cc}a&b\\c&d\end{array}\right)=\frac{(1+a^2)d-abc}{\sqrt{T\left(\begin{array}{cc}a&b\\c&d\end{array}\right)}}.\end{split}\ee

Then these functions $A_i^j,i,j=x,y$ are $C^\infty$ near the origin, and for any compact neighborhood $K$ near the origin, all its derivatives are uniformly controlled by some constant depending on $K$.

The system (4.24) becomes
\be \left\{\begin{split}D_x(A_x^x(\nabla f))+D_y(A_y^x(\nabla f))=0,\\
D_x(A_x^y(\nabla f))+D_y(A_y^y(\nabla f))=0.\end{split}\right.\ee

We differentiate (4.26) with respect to $x$, we have
\be\left\{\begin{split}D_x&[D_aA_x^x(\nabla f)\cdot D_x u_x+D_bA_x^x(\nabla f)\cdot D_xv_x+D_cA_x^x(\nabla f)\cdot D_yu_x+D_d A_x^x(\nabla f)\cdot D_yv_x]+\\
&D_y[D_aA_y^x(\nabla f)\cdot D_x u_x+D_bA_y^x(\nabla f)\cdot D_x v_x+D_cA_y^x(\nabla f)\cdot D_y u_x+D_d A_y^x(\nabla f)\cdot D_y v_x]=0,\\
D_x&[D_aA_x^y(\nabla f)\cdot D_xu_x+D_bA_x^y(\nabla f)\cdot D_xv_{x}+D_cA_x^y(\nabla f)\cdot D_yu_{x}+D_d A_x^y(\nabla f)\cdot D_yv_{x}]+\\
&D_y[D_aA_y^y(\nabla f)\cdot D_x u_x+D_bA_y^y(\nabla f)\cdot D_x v_x+D_cA_y^y(\nabla f)\cdot D_y u_x+D_d A_y^y(\nabla f)\cdot D_y v_x]=0.
\end{split}\right.\ee
This means that the function $(u_x,v_x)$ satisfies the above system, with coefficient matrix
\be A(\nabla f)=\left(\begin{array}{cccc}
D_aA_x^x(\nabla f)&D_cA_x^x(\nabla f)&D_aA_x^y(\nabla f)&D_cA_x^y(\nabla f)\\
D_aA_y^x(\nabla f)&D_cA_y^x(\nabla f)&D_aA_y^y(\nabla f)&D_cA_y^y(\nabla f)\\
D_bA_x^x(\nabla f)&D_d A_x^x(\nabla f)&D_bA_x^y(\nabla f)&D_d A_x^y(\nabla f)\\
D_bA_y^x(\nabla f)&D_dA_y^x(\nabla f)&D_bA_y^y(\nabla f)&D_d A_y^y(\nabla f)
\end{array}\right).\ee 

We calculate the partial derivates of $A_i^j,i,j=x,y$, for $\left(\begin{array}{cc}a&b\\c&d\end{array}\right)\in M_2(\R)$, and get 
\be A\left(\begin{array}{cc}a&b\\c&d\end{array}\right)=\left(\begin{array}{cccc}
\frac{1+d^2-(A_x^x)^2}{\sqrt T}&\frac{-A_x^xA_y^x-bd}{\sqrt T}&\frac{-A_x^xA_x^y-cd}{\sqrt T}&\frac{-A_x^yA_y^x+2bc-ad}{\sqrt T}\\
\frac{-A_y^xA_x^x-bd}{\sqrt T}&\frac{1+b^2-(A_y^x)^2}{\sqrt T}&\frac{-A_y^yA_x^x+2ad-bc}{\sqrt T}&\frac{-A_y^yA_y^x-ab}{\sqrt T}\\
\frac{-A_x^xA_x^y-cd}{\sqrt T}&\frac{-A_x^xA_y^y+2ad-bc}{\sqrt T}&\frac{1+c^2-(A_x^y)^2}{\sqrt T}&\frac{-A_x^yA_y^y-ac}{\sqrt T}\\
\frac{-A_y^xA_x^y+2bc-ad}{\sqrt T}&\frac{-A_y^xA_y^y-ab}{\sqrt T}&\frac{-A_y^yA_x^y-ac}{\sqrt T}&\frac{1+a^2-(A_y^y)^2}{\sqrt T}
\end{array}\right).\ee

We can observe that when $a,b,c,d$ are small enough, $A\left(\begin{array}{cc}a&b\\c&d\end{array}\right)$ satisfies the strong elliptic condition (3.12) in \cite{Gia}, hence the coefficient matrix $A(\nabla f)$ of (4.29) satisfies the strong elliptic condition, when $\mu$ is small. Moreover the $C^{0,\sigma}$ norm of $A(\nabla f)$ is also controlled by $||f||_{C^{1,\sigma}}$, and hence by $\mu$. 

Hence for the function $(u_x,v_x)$, by Caccioppoli's inequality (cf.\cite{Gia} Theorem 4.4), we have
\be ||\nabla(u_x,v_x)||_{L^2(\overline B(0,\frac 78))}\le C||(u_x,v_x)||_{L^2(\overline B(0,\frac 89))}\le C||f||_{C^{1,\sigma}},\ee
where $C$ depends on the $C^{0,\sigma}$ norm of the coefficient matrix $A(\nabla f)$, hence by $||f||_{C^{1,\sigma}}$, hence by $\mu$.

Then by the Schauder estimates (Theorem 5.17 of \cite{Gia}), we have
\be ||\nabla(u_x,v_x)||_{C^{0,\sigma}(\overline B(0,\frac 67))}\le C(\mu)||\nabla(u_x,v_x)||_{L^2(\overline B(0,\frac 78))}\le C||f||_{C^{1,\sigma}}\le C_2'(\mu),\ee
where $C_2'(\mu)\to 0$ while $\mu\to 0$.

We differentiate the system (4.26) with respect to y, and get the same estimation
\be ||\nabla(u_y,v_y)||_{C^{0,\sigma}(\overline B(0,\frac 67))}\le C_2'(\mu).\ee

Hence we get
\be ||f||_{C^{2,\sigma}(\overline B(0,\frac 67))}\le C_2(\mu),\ee
with $\lim_{\mu\to 0}C_2(\mu)=0$. 

We still need to estimate $\nabla^3 f$. For this we differentiate the system (4.27). We set $g_1=u_x, g_2=v_x$, and for $i=x,y, j=1,2$, set $p_{x1}=a, p_{x2}=b, p_{y1}=c, p_{y2}=d$. Then (4.27) becomes,
\be \sum_{\a=x,y}D_\a(\sum_{i=x,y,j=1,2} D_{P_{ij}} A_\a^\beta(\nabla f)\cdot D_i g_j)=0, \mbox{ for }\beta=x,y.\ee

Now we differentiate it with respect to $s$, for $s\in\{x,y\}$, and get
\be \sum_{\a=x,y}D_\a(\sum_{i=x,y,j=1,2} D_{P_{ij}} A_\a^\beta(\nabla f)\cdot D_i (D_sg_j))+\sum_{\a=x,y}D_\a(\sum_{i=x,y,j=1,2} D_{P_{ij}} D_sA_\a^\beta(\nabla f)\cdot D_i g_j)=0,\ee
$\beta=x,y$. I.e. the function $(D_s g_1,D_s g_2)$ satisfies the elliptic system
\be \sum_{\a=x,y}D_\a(\sum_{i=x,y,j=1,2} D_{P_{ij}} A_\a^\beta(\nabla f)\cdot D_i (D_sg_j))=-\sum_{\a=x,y}D_\a(\sum_{i=x,y,j=1,2} D_{P_{ij}} D_sA_\a^\beta(\nabla f)\cdot D_i g_j).\ee

Notice that the left hand side of the system is exactly the same as (4.34), hence the function $(D_sg_1,D_sg_2)$ is a solution to the elliptic system
\be  \sum_{\a=x,y}D_\a(\sum_{i=x,y,j=1,2} D_{P_{ij}} A_\a^\beta(\nabla f)\cdot D_i (D_sg_j))=-\sum_{\a=x,y}D_\a(\sum_{i=x,y,j=1,2} B_{i,j}^{\a,\beta}),\ee
where $B_{i,j}^{\a,\beta}=D_{P_{ij}} D_sA_\a^\beta(\nabla f)\cdot D_i g_j$, hence $||B_{i,j}^{\a,\beta}||_{C^{0,\sigma}}$ is controlled by $||f||_{C^{2,\sigma}}$, which is controlled by $C_2(\mu)$, and is small.

We apply again the Caccioppoli's inequality for $(D_sg_1,D_sg_2)$, and get
\be \begin{split}||\nabla(D_sg_1,D_sg_2)||_{L^2(\overline B(0,\frac 56))}&\le C(||(D_sg_1,D_sg_2)||^2_{L^2(\overline B(0,\frac 78)}+||\sum_{\a=x,y,\beta=x,y}\sum_{i=x,y,j=1,2} B_{i,j}^{\a,\beta}||^2_{L^2(\overline B(0,\frac 67)})^\frac12\\
&\le C(||\nabla f||^2_{L^2(\overline B(0,\frac 78)})\le C_3'(\mu),\end{split}\ee
with $\lim_{\mu\to 0}C_3'(\mu)=0$.

Then we apply again the Schauder estimates (Theorem 5.17 of \cite{Gia}), and get
\be \begin{split}||\nabla (D_sg_1,D_sg_2)||_{C^{0,\sigma}(\overline B(0,\frac 45))}&\le C(||\nabla (D_sg_1,D_sg_2)||_{L^2(\overline B(0,\frac 56))}+||\sum_{i=x,y,j=1,2} B_{i,j}^{\a,\beta}||_{C^{0,\sigma}(\overline B(0,\frac 67))})\\
&\le C_3''(\mu), \mbox { for }s=x,y,
\end{split}\ee
with $\lim_{\mu\to 0}C_3''(\mu)=0$.

Recall that $(g_1,g_2)=(u_x,v_x)$. We repeat the same argument for $(u_y,v_y)$, and altogether we have
\be ||\nabla^3 f||_{C^{0,\sigma}(\overline B(0,\frac 45))}\le C_3(\mu),\ee
with  $\lim_{\mu\to 0}C_3(\mu)=0$.

Combining (4.1), (4.15), (4.33) and (4.40), we have that for any $\mu$ small, there exists a constant $C(\mu)$, with $\lim_{\mu\to 0}C(\mu)=0$, such that if $f$ is a minimal graph on $B(0,1)$, with 
\be\max\{||f|_{\partial B(0,1)}||L^\infty, ||Df|_{\partial B(0,1)}||L^\infty\}\le\mu,\ee then
\be \max_{0\le i\le 3}||\nabla^i f||_{L^\infty(B(0,\frac 34)}\le C(\mu).\ee

Thus we complete the proof of Proposition 4.2.\qed

\section{Estimates for perturbations around a minimal graph}

Denote by $B=B(0,1)\cap \R^2$ the unit disc in $\R^2$. Let $q\in B(0,\frac{1}{100})$, and set $B_r=B(q,r)$ for $r>0$. Fix any $\e$ and $l$ less than $10^{-4}$, let $\mu<10^{-4}$ be small. (Here in this section the three are independent; in the next section, $l$ will be chosen first, and then $\e$ will  depend on $l$, and both will be fixed at the beginning, while $\mu$ will be supposed to be much smaller than these two, and will be decided later.) Let $f$ be a function from $\overline B$ to $\R^2$ whose graph $\Sigma_f=\{(x,f(x));x\in \overline B\}\subset\R^4$ is a minimal submanifold in $\R^4$, with $||f|_{\partial B}||_{C^1}\le\mu$. Let $h$ be a $C^1$ function from $A_r:=\overline B\bs B_r$ to $\R^2$ with $h|_{\partial B}=0$, Lip $h\le l$, and there exists a vector $M\in\R^2$ such that for any $x\in\partial B_r$, $|h(x)-M|\le\e r$. Denote by $\Sigma_{f+h}$ the graph of $f+h$ on the annulus $A_r$.

\begin{pro}Take all the notations and assumptions above, then
\be  H^2(\Sigma_{f+h})-H^2(\Sigma_f)\ge\frac14\int_{A_r}|\nabla h|^2-C r^2(\mu+\mu\e+C_0(\mu)),\ee
where $\lim_{\mu\to 0}C_0(\mu)=0.$
\end{pro}

\nd Now let us compare $\Sigma_{f+h}$ and $\Sigma_f$ above $A_r$. We have
\be\begin{split}H^2(\Sigma_{f+h})&-H^2(\Sigma_f)=\int_{A_r}\sqrt{1+S(f+h)}-\sqrt{1+S(f)}\\
&=\int_{A_r}\sqrt{1+S(f)}(\sqrt{\frac{1+S(f+h)}{1+S(f)}}-1)\\
&=\int_{A_r}\sqrt{1+S(f)}(\sqrt{1+\frac{S(f+h)-S(f)}{1+S(f)}}-1).
\end{split}\ee
But
\be\begin{split}&S(f+h)-S(f)=[|\nabla(f+h)|^2-|\nabla f|^2]+[(\det\nabla (f+h))^2-(\det\nabla f)^2]\\
&=[2<\nabla f,\nabla h>+|\nabla h|^2]+[<(\nabla f)^*,\nabla h>+\det\nabla h)][2\det\nabla f+\det\nabla h+<(\nabla f)^*,\nabla h>].\end{split}\ee
Notice that $|\nabla f|<2\mu$, $|(\nabla f)^*|<2\mu$ is small, and $|\det\nabla f|\le|\nabla f|^2$, $|\det\nabla h|\le|\nabla h|^2$,
therefore $|S(f+h)-S(f)|<1$ since $|\nabla h|<l$ is small. But $S(f)>0$, hence $|\frac{S(f+h)-S(f)}{1+S(f)}|<1$. For any $|x|<1$ we have 
\be 1+x=(1+\frac x2)^2-\frac{x^2}{4}\ge(1+\frac x2-\frac{x^2}{4})^2,\ee
hence
\be \sqrt{1+\frac{S(f+h)-S(f)}{1+S(f)}}\ge 1+\frac12 \frac{S(f+h)-S(f)}{1+S(f)}-\frac14(\frac{S(f+h)-S(f)}{1+S(f)})^2,\ee
which gives
\be \begin{split}H^2(\Sigma_{f+h})-H^2(\Sigma_f)&
\ge\int_{A_r}\sqrt{1+S(f)}(\frac12 \frac{S(f+h)-S(f)}{1+S(f)}-\frac14(\frac{S(f+h)-S(f)}{1+S(f)})^2)\\
&=\frac12\int_{A_r}\frac{S(f+h)-S(f)}{\sqrt{1+S(f)}}-\frac14\int_{A_r}\frac{(S(f+h)-S(f))^2}{(1+S(f))^\frac32}
.\end{split}\ee
For the first term, by (5.4),
\be\begin{split} \frac12\int_{A_r}&\frac{S(f+h)-S(f)}{\sqrt{1+S(f)}}=\frac12\int_{A_r}\frac{2<\nabla f,\nabla h>+|\nabla h|^2+2\det\nabla f<(\nabla f)^*,\nabla h>}{\sqrt{1+S(f)}}+\\
&\frac12\int_{A_r}\frac{2<(\nabla f)^*,\nabla h>\det\nabla h+<(\nabla f)^*,\nabla h>^2+2\det\nabla h\det\nabla f+|\det\nabla h|^2}{\sqrt{1+S(f)}}\\
&\ge\int_{A_r}\frac{<\nabla f,\nabla h>+\frac 12|\nabla h|^2+\det\nabla f<(\nabla f)^*,\nabla h>}{\sqrt{1+S(f)}}-(2\mu+l^2)\int_{A_r}|\nabla h|^2
\end{split}\ee
But $S(f)\le 5\mu^2$, hence $\frac{1}{1+S(f)}\ge \frac89$, hence we have
\be\frac12\int_{A_r}\frac{S(f+h)-S(f)}{\sqrt{1+S(f)}}\ge\int_{A_r}<\frac{\nabla f+\det\nabla f(\nabla f)^*}{\sqrt{1+S(f)}},\nabla h>+\frac13\int_{A_r}|\nabla h|^2.\ee
By (4.23), and the hypothesis that $h|_{\partial B}=0$, we have
\be \begin{split}&\int_{A_r}<\frac{\nabla f+\det\nabla f(\nabla f)^*}{\sqrt{1+S(f)}},\nabla h>\\
&=\int_{\partial A_r}<h,[\vec n\cdot\frac{\nabla f+\det\nabla f(\nabla f)^*}{\sqrt{1+S(f)}}]>-\int_{A_r}<div(\frac{\nabla f+\det\nabla f(\nabla f)^*}{\sqrt{1+S(f)}}),h>\\
&=-\int_{\partial B_r}<h,[\vec n\cdot\frac{\nabla f+\det\nabla f(\nabla f)^*}{\sqrt{1+S(f)}}]>\\
&-\int_{\partial B_r}<(M+h-M),[\vec n\cdot\frac{\nabla f+\det\nabla f(\nabla f)^*}{\sqrt{1+S(f)}}]>\\
&=-<M,\int_{\partial B_r}[\vec n\cdot\frac{\nabla f+\det\nabla f(\nabla f)^*}{\sqrt{1+S(f)}}]>+\int_{\partial B_r}<(M-h),[\vec n\cdot\frac{\nabla f+\det\nabla f(\nabla f)^*}{\sqrt{1+S(f)}}]>.
\end{split}\ee

For the second term of (5.10), since $|M-h|\le\e r, $ Lip$f\le\mu$, and $|det\nabla f|\le 2|\nabla f|^2\le 2\mu^2\le \mu$ since $\mu$ is small, we have
\be |\int_{\partial B_r}<(M-h),[\vec n\cdot\frac{\nabla f+\det\nabla f(\nabla f)^*}{\sqrt{1+S(f)}}]>|
\le \int_{\partial B_r}\e r (2\mu)\le 4\pi\mu\e r^2.\ee

For the first term of (5.10), first by Taylor expansion at the point $0$, we have, for any $x\in\partial B_r$,
\be \nabla f(x)=\nabla f(0)+x\cdot\nabla^2 f(0)+o_1(r), \ee
\be ( \nabla f)^*(x)=(\nabla f)^*(0)+x\cdot\nabla(\nabla)^*f(0)+o_2(r),\ee
\be \det(\nabla f)(x)=\det(\nabla f)(0)+x\cdot \nabla\det(\nabla f)(0)+o_3(r),\ee
\be  \frac{1}{\sqrt{1+S(f)}}(x)=\frac{1}{\sqrt{1+S(f)}}(0)+x\cdot \nabla(\frac{1}{\sqrt{1+S(f)}})(0)+o_4(r)\ee
where $|o_1(r)|\le r^2||\nabla^3 f||_{L^\infty (B(0,r))} $, $|o_2(r)|\le r^2||\nabla^3 f||_{L^\infty( B(0,r))}$, $|o_3(r)|\le r^2||\nabla^2 det(\nabla f)||_{L^\infty (B(0,r))}$, $|o_4(r)|\le  r^2||\nabla^2 (\frac{1}{\sqrt{1+S(f)}})||_{L^\infty (B(0,r))}$.

Hence we have
\be \begin{split}&\frac{\nabla f+\det\nabla f(\nabla f)^*}{\sqrt{1+S(f)}}=\\
&\{\nabla f(0)+x\cdot\nabla^2 f(0)+o_1(r)+[\det(\nabla f)(0)+x\cdot \nabla\det(\nabla f)(0)+o_3(r)][(\nabla f)^*(0)+x\cdot\nabla(\nabla)^*f(0)+o_2(r)]\}\\
&[\frac{1}{\sqrt{1+S(f)}}(0)+x\cdot \nabla(\frac{1}{\sqrt{1+S(f)}})(0)+o_4(r)]\\
&=\{[\nabla f(0)+\det(\nabla f)(0)(\nabla f)^*(0)]+x\cdot[\nabla^2 f(0)+\nabla\det(\nabla f)(0)(\nabla f)^*(0)+\det(\nabla f)(0)\nabla(\nabla)^*f(0)]+o(r)\}\\
&[\frac{1}{\sqrt{1+S(f)}}(0)+x\cdot \nabla(\frac{1}{\sqrt{1+S(f)}})(0)+o(r)]\\
&=[\nabla f(0)+\det(\nabla f)(0)(\nabla f)^*(0)]\frac{1}{\sqrt{1+S(f)}}(0)\\
&+x\cdot \frac{1}{\sqrt{1+S(f)}}(0)[\nabla^2 f(0)+\nabla\det(\nabla f)(0)(\nabla f)^*(0)+\det(\nabla f)(0)\nabla(\nabla)^*f(0)]\\
&+[\nabla f(0)+\det(\nabla f)(0)(\nabla f)^*(0)][x\cdot\nabla(\frac{1}{\sqrt{1+S(f)}})(0)]+o(r),\end{split}.\ee
where all the $o(r)$ in (5.16) satisfied that $|o(r)|\le C_0r^2$, where 
\be C_0=C(||\nabla f||_{L^\infty B(0,r)},||\nabla^2 f||_{L^\infty B(0,r)},||\nabla^3 f||_{L^\infty B(0,r)})\ee
 tends to 0 as $||\nabla f||_{L^\infty B(0,r)},||\nabla^2 f||_{L^\infty B(0,r)},||\nabla^3 f||_{L^\infty B(0,r)}$ tend to 0.

Therefore,
\be \begin{split}&|-<M,\int_{\partial B_r}[\vec n\cdot\frac{\nabla f+\det\nabla f(\nabla f)^*}{\sqrt{1+S(f)}}]>|\\
&\le |<M,\int_{\partial B_r}[\vec n\cdot[\nabla f(0)+\det(\nabla f)(0)(\nabla f)^*(0)]\frac{1}{\sqrt{1+S(f)}}(0)]>|\\
&+|<M,\int_{\partial B_r}[\vec n\cdot(x\cdot \frac{1}{\sqrt{1+S(f)}}(0)[\nabla^2 f(0)+\nabla\det(\nabla f)(0)(\nabla f)^*(0)+\det(\nabla f)(0)\nabla(\nabla)^*f(0)])]|\\
&+|<M,\int_{\partial B_r}\{\vec n\cdot[\nabla f(0)+\det(\nabla f)(0)(\nabla f)^*(0)][x\cdot\nabla(\frac{1}{\sqrt{1+S(f)}})(0)]\}>|+|<M,\int_{\partial B_r}o(r)>|.
\end{split}\ee
For the first term of (5.18), since $[\nabla f(0)+\det(\nabla f)(0)(\nabla f)^*(0)]\frac{1}{\sqrt{1+S(f)}}(0)$ is a constant matrix, which we denote by $V$, and hence we have
\be <M,\int_{\partial B_r}\vec n\cdot[\nabla f(0)+\det(\nabla f)(0)(\nabla f)^*(0)]\frac{1}{\sqrt{1+S(f)}}(0)>=<M,(\int_{\partial B_r} \vec n)\cdot V>=0\ee
because $\int_{\partial B_r} \vec n=0$. 

For the second and third term of (5.16), notice that $|x|=r$, $\nabla f\le \mu$, hence their sum is less than
\be C\mu r^2+C|\nabla^2 f(0)| r^2\le (C\mu+CC_0)r^2,\ee
where $C_0$ is as in (5.17) and $C$ does not depend on $\mu,\e$.

For the last, by the previous control on $o(r)$, this term is less than $ C_0r^3$.

Altogether we have
\be |-<M,\int_{\partial B_r}[\vec n\cdot\frac{\nabla f+\det\nabla f(\nabla f)^*}{\sqrt{1+S(f)}}]>|\le C r^2(\mu+\ C_0).\ee

Combining with (5.11) and (5.9), we have
\be \frac12\int_{A_r}\frac{S(f+h)-S(f)}{\sqrt{1+S(f)}}\ge \frac13\int_{A_r}|\nabla h|^2-C r^2(\mu+\mu\e+\ C_0),\ee
where $C$ does not depend on $\mu,l$ and $\e$.

Recall that this is the estimation for the first term of the last line in (5.7). Now we treat its second term.

By (5.4), we have
\be\begin{split} &|S(f+h)-S(f)|\\
=&|[2<\nabla f,\nabla h>+|\nabla h|^2]+[<(\nabla f)^*,\nabla h>+\det\nabla h)][2\det\nabla f+\det\nabla h+<(\nabla f)^*,\nabla h>]|\\
\le & 2|\nabla f||\nabla h|+|\nabla h|^2+(|(\nabla f)^*||\nabla h|+|\nabla h|^2|][2|\nabla f|^2+|\nabla h|^2+|(\nabla f)^*||\nabla h|]\\
\le &C(|\nabla f||\nabla h|+|\nabla h|^2)\le C\mu|\nabla h|+C|\nabla h|^2,\end{split}\ee
therefore the second term of (5.7) verifies
\be\begin{split}&-\frac14\int_{A_r}\frac{(S(f+h)-S(f))^2}{(1+S(f))^\frac32}\ge -\frac14\int_{A_r}(S(f+h)-S(f))^2\\
&\ge-\frac14\int_{A_r} (C\mu|\nabla h|+C|\nabla h|^2)\ge -C(\mu^2+||\nabla h||_\infty^2)\int_{A_r}|\nabla h|^2.\end{split}\ee

On combining (5.7), (5.22) and (5.24) we get
\be \begin{split}H^2(\Sigma_{f+h})-H^2(\Sigma_f)&\ge \frac13\int_{A_r}|\nabla h|^2-C r^2(\mu+\mu\e+\ C_0)-C(\mu^2+||\nabla h||_\infty^2)\int_{A_r}|\nabla h|^2\\
&\ge (\frac 13-C\mu^2-Cl^2)\int_{A_r}|\nabla h|^2-C r^2(\mu+\mu\e+\ C_0).
\end{split}\ee
But Lip $h<l$ is small, hence we have
\be H^2(\Sigma_{f+h})-H^2(\Sigma_f)\ge\frac14\int_{A_r}|\nabla h|^2-C r^2(\mu+\mu\e+C_0).\ee

Now we apply Proposition 4.2, and get that when $r<\frac 34$ and $\mu$ is small enough, 
\be C_0=C_0(||\nabla f||_{L^\infty B(0,r)},||\nabla^2 f||_{L^\infty B(0,r)},||\nabla^3 f||_{L^\infty B(0,r)})=C_0(C(\mu))=C_0(\mu),\ee with $\lim_{\mu\to 0}C_0(\mu)=0$. Thus we have
\be  H^2(\Sigma_{f+h})-H^2(\Sigma_f)\ge\frac14\int_{A_r}|\nabla h|^2-C r^2(\mu+\mu\e+C_0(\mu)).\ee\qed

\section{Conclusion}

Now return to our set $E$. Recall that $\a$ is a pair of angles larger than $\theta_1'>\frac\pi 3$. $E$ is a reduced closed set that is minimal in $B(0,1)$, which contains no $2\P$ type point in $B(0,\frac{1}{100})$. 

Set $l=10^{-3}$, and suppose that $d_{0,1}^\a<\mu<\min\{\frac{\e_0}{10},\frac l2\}$, $\mu$ is to be decided later. 
%
%

We apply Proposition 2.11 to $E$, with $\e'=\min\{\e_{\frac l2}, 10^{-4}\}$, (where $\e_{\frac l2}$ corresponds to $\frac l2$ in Proposition 2.20), and get our $o_E$ and $r_E$. Then $r_E<\frac14$.

Let $\gamma^i,g^i,$ as in Section 3. Suppose that 
\be ||\gamma^i||_{C^1}\le \mu, i=1,2.\ee

By Theorem 4.1 and 4.2 of \cite{LaOs}, for each $i$ there exists a function $f^i:\overline B(0,\frac 12)\cap P_\a^i\to{P_\a^i}^\perp$, whose graphs $\Sigma^i=\Sigma_{f^i}=\{(x,f(x)):x\in \overline B(0,\frac 12)\cap P_\a^i\}\subset\R^4$ are minimal surfaces. Denote by $B^i(x,r)=B(x,r)\cap P_\a^i$.

On the other hand, we want to show the part of $E$ in the annulus $D_\a(o_E,r_E)\bs D_\a(o_E,\frac14 r_E)$ is far from any translation of $P_\a$. Recall that Proposition 2.11 says that $E$ is $\e' r_E$ far from any translation of $P_\a$ in the ball $D_\a(o_E,r_E)$. So for having a relatively big distance in the annulus, we simply use a compactness argument, and can get the following proposition. (See \cite{2p} for the proof).

\begin{pro}[cf.\cite{2p}, Corollary 8.24]
For every $\e>0$, there exists $0<\d<\e$, and $0<\theta_0<\frac\pi2$, which do not depend on $\e$, with the following properties. If $\theta_0<\theta<\frac\pi2$, and if $E$ is minimal in $D_\theta(0,1)$ and is $\delta$ near $P_\theta$ in $D_\theta(0,1)\bs D_\theta(0,\frac 14)$, and moreover
\be p_\theta^i(E)\supset P_\theta^i\cap B(0,\frac 34),\ee
then $E$ is $\e$ near $P_\theta$ in $D_\theta(0,1)$.
\end{pro}

Let $\d'$ be the $\d$ corresponding to $\e'$ in Proposition 6.2, we know that $E$ is not $\d' r_E$ near any translation of $P_\a$ in $D_\a(o_E,r_E)\bs D_\a (o_E,\frac 14 r_E)$. On the other hand, by definition of $o_E$ and $r_E$, we know that the $\e'$-process does not stop at the scale $2r_E$, thus by Proposition 2.20, $E\cap D_\a(o_E,r_E)\bs D_\a (o_E,\frac 14 r_E)$ is composed of two fine $C^1$ graphs $G^1,G^2$ of two functions $g^i, i=1,2$ on $P_\a^i\cap D_\a(o_E,r_E)\bs D_\a (o_E,\frac 14 r_E)$ respectively. Thus $G^1\cup G^2$ is not $\d'r_E$ near any translation of $P_\a$, there exists $i=1,2$ such that $G^i$ is not $\d'$ near any translation of $P_\a^i$ in $D_\a(o_E,r_E)\bs D_\a (o_E,\frac 14 r_E)$. Suppose this is the case for $i=1$. 

Denote by $g=g^1, f=f^1,$ and $h=g-f$. We want to apply Proposition 5.1 to $f$ and $h$, with $B(q,r)=B^1(o_E,\frac 14r_E)$ (hence $q=o_E,r=\frac 14 r_E$). Recall that we have set $\e'\le \e_{\frac l2}$, hence $|\nabla g|$ is smaller than $\frac l2$, which gives $|\nabla h|=|\nabla (g-f)|$ is smaller than $|\nabla g|+|\nabla f|<\frac l2+\mu<l$ cause $\mu$ is supposed to be less than $\frac l2$.

 Also, by Proposition 2.11, $G^1$ is still $2\e' r_E$ near some translation of $P_\a^1$, hence there exists $M_g\in {P_\a^1}^\perp$ such that $|g(x)-M_g|\le 2\e' r_E=8\e' r$. But $f$ is $\mu$-Lipschitz, hence there exists $M_f$ such that $|f(x)-M_f|\le C\mu r$ on $\partial B(q,r)$, which gives $|h-(M_g+M_f)|\le 9\e' r <10^{-3}r$ on $\partial B(q,r)$, when $\mu$ is small. 
  
Now we can apply Proposition 5.1, and get
\be H^2(G^1)-H^2(\Sigma^1\bs C^1(o_E,\frac 14 r_E))=H^2(\Sigma_{f+h})-H^2(\Sigma_f)\ge \frac14\int_{A_r}|\nabla h|^2-C r^2(\mu+\e'\mu+C_0(\mu)),\ee
with $A_r=B^1(0,\frac 12)\bs B(q,r)$.

Now we want to estimate $\int_{A_r}|\nabla h|^2$. Recall that on $B^1(o_E,r_E)\bs B^1(o_E,\frac 14 r_E)$, the graph of $g$ is $\d' r_E$ far from any translation of $P_\a^1$. On the other hand $f$ is $\mu$-Lipschitz, hence when $\mu$ is small, the graph of $h=g-f$ is $\frac12 \d' r_E$ far from any translation of $P_\a^1$.

Firstly we cite here two lemmas for estimating the Dirichlet's energy of our perturbation function $h$.

\begin{lem}[cf.\cite{2p}, Corollary 7.23]\label{cas1}Let $r_0>0$, $q\in\R^2$ be such that $r_0<\frac12d(q,\partial B(0,1))$, suppose $u_0\in C^1(\partial B(q,r_0)\cap \R^2,\R)$, and denote by $m(u_0)=\frac{1}{2\pi r_0}\int_{\partial B(q,r_0)}u_0$ its average. 

Then for all $u\in C^1((\overline{B(0,1)}\backslash B(q,r_0))\cap \R^2, \R)$ that satisfies
\be  u|_{\partial B(q,r_0)}=u_0\ee 
we have \be \int_{B(0,1)\backslash B(q,r_0)}|\nabla u|^2\ge \frac14r_0^{-1}\int_{\partial B(q,r_0)}|u_0-m(u_0)|^2.\ee 
\end{lem}

\begin{lem}[cf.\cite{2p}, Corollary 7.36]\label{level}
For all $0<\epsilon<1$, there exists $C=C(\epsilon)>100$ such that if $0<r_0<1,$ $u\in C^1(\ B(0,1)\backslash B(0,r_0),\R)$ and
\be u|_{\partial B(0,r_0)}>\delta r_0-\frac{\delta r_0}{C}\ and\ \ u|_{\partial B(0,1)}<\frac{\delta r_0}{C}\ee 
then \be \int_{B(0,1)\backslash B(0,r_0)} |\nabla u|^2\ge\epsilon \frac{2\pi\delta^2r_0^2}{|\log r_0|}.\ee 
\end{lem}

Then denote by $P=P_\a^1$ for short. Denote by $D=D_\a$. Then $h$ is a map from $P$ to $P^\perp$, and is therefore from $\R^2$ to $\R^2$. Write $h=(\varphi_1,\varphi_2)$, where $\varphi_i:\R^2\to \R$. Then since the graph of $h$ is $\frac12 \d' r_E$ far from all translation of $P$, there exists $j\in\{1,2\}$ such that
\be \sup_{x,y\in P\cap D(o_E, r_E)\backslash D(o_E, \frac 14 r_E)}|\varphi_j(x)-\varphi_j(y)|\ge \frac 14 r_E\d'.\ee

Suppose this is true for $j=1$. Denote by 
\be K=\{(z, \varphi_1(z)): z\in (D(0,\frac 12)\backslash D(o_k,\frac14 r_E))\cap P\},\ee 
then
\be \begin{split}K\mbox{ is the orthogonal }&\mbox{projection of }G^1\cap D(0,\frac 12)\\
&\mbox{ on a 3-dimensional subspace of }\R^4.\end{split}\ee 


For $\frac 14r_E\le s\le r_E$, define 
\be \Gamma_s=K\cap p^{-1}(\partial D(o_E, s)\cap P)=\{(x,\varphi_1(x))|x\in\partial D(o_E, s)\cap P\}\ee  the graph of $\varphi_1$ on $\partial D(o_E, s)\cap P$. 

We know that the graph of $\varphi_1$ is $\frac14\d'r_E$ far from $P$ in $D(o_E,r_E)\bs D(o_E,\frac 14r_E)$; then there are two cases:

1st case: there exists $t\in[\frac 14r_E, r_E]$ such that
\be \sup_{x,y\in \Gamma_t}\{|\varphi_1(x)-\varphi_1(y)|\}\ge \frac{\d'}{C}r_E,\ee 
where $C=4C(\frac 12)$ is the constant of Lemma \ref{level}.

Then there exists $a,b\in \Gamma_t$ such that $|\varphi_1(a)-\varphi_1(b)|>\frac{\d'}{ C} r_E\ge\frac{\d'}{ C}t$. Since $||\nabla \varphi_1||_\infty\le ||\nabla \varphi||_\infty<1$, we have 
\be \int_{\Gamma_t}|\varphi_1-m(\varphi_1)|^2\ge \frac{t^3\d'^3}{4C^3}=(\frac43 t\d')^3(\frac{27}{4^4C^3}).\ee 
Now in $D(0,\frac 12)$ we have $d(0,o_E)<6\e'\le 10\e' \cdot\frac12$, and $s<r_E<\frac 18<\frac 12\times\frac 12$, therefore we can apply Lemma \ref{cas1} and obtain
\be \int_{(D(0,\frac 12)\backslash D(o_E, t))\cap P}|\nabla \varphi_1|^2\ge C(\d')t^2\ge C_1(\d')r_E^2.\ee 

2nd case: for all $\frac 14r_E\le s\le r_E$, 
\be \sup_{x,y\in \Gamma_s}\{|\varphi_1(x)-\varphi_1(y)|\}\le \frac{\d'}{C}r_E.\ee 
However, since
\be\begin{split}
\frac12 r_E\d'&\le\sup\{|\varphi_1(x)-\varphi_2(y)|:x,y\in P\cap D(o_E,r_E)\backslash D(o_E,\frac14 r_E)\}\\
&=\sup\{|\varphi_1(x)-\varphi_2(y)|:s,s'\in[\frac14r_E,r_E],x\in\Gamma_s,y\in\Gamma_{s'}\},
\end{split}\ee
there exist $\frac 14r_E\le t<t'\le r_E$ such that
\be \sup_{x\in \Gamma_{t},y\in\Gamma_{t'}}\{|\varphi_1(x)-\varphi_1(y)|\}\ge \frac12r_E\d'. \ee 
Fix $t$ and $t'$, and without loss of generality, suppose that
\be \sup_{x\in \Gamma_{t},y\in\Gamma_{t'}}\{\varphi_1(x)-\varphi_1(y)\}\ge \frac14r_E\d' .\ee 
Then
\be \inf_{x\in\Gamma_{t}}\varphi_1(x)-\sup_{x\in\Gamma_{t'}}\varphi_1(x)\ge \frac14r_E\d'-2\frac{\d'}{C}r_E=(1-\frac{2}{C(\frac12)})\frac{\d'}{4}r_E\ge (1-\frac{2}{C(\frac12)})\frac{\d'}{2}t'\ee 
because $C=4 C(\frac 12)$.

Now look at what happens in the ball $D(o_E,t')\cap P$. Apply Lemma \ref{level} to the scale $t'$, we get
\be \int_{(D(o_E,t')\backslash D(o_E, t))\cap P}|\nabla \varphi_1|^2\ge C(\d',\frac 12)\frac{\pi(\frac{\d'}{2})^2 t'^2}{\log\frac{t'}{t}}.\ee 
Then since $\frac{t'}{t}\le 4, t'>t>\frac 14r_E$, we have 
\be \int_{((D(o_E,t')\backslash D(o_E, t))\cap P}|\nabla \varphi_1|^2\ge C_2(\d')r_E^2.\ee 

So in both cases, there exists a constant $C=C_5(\d')=\min\{C_1(\d'), C_2(\d')\}$, which depends only on $\d'$, such that
\be \int_{(D(0,\frac 12)\backslash D(o_E, t_E))\cap P}|\nabla \varphi_1|^2\ge C_5(\d')r_E^2.\ee 

On the other hand, since $|\nabla \varphi_1|\le|\nabla h|<1$, we have
\be \int_{A_r}|\nabla h|^2=\int_{(D(0,\frac 12)\backslash D(o_E, t_E))\cap P}|\nabla h|^2\ge C_5(\d')r_E^2.\ee

Thus by (6.4),
\be H^2(G^1)-H^2(\Sigma^1\bs C^1(o_E,\frac 14 r_E))\ge C_5(\d')r_E^2-Cr_E^2(\mu+\e'\mu+C_0(\mu)).\ee

We apply also Proposition 5.1 to $i=2$, where all the verifications for $g^2$, $f^2$, $h^2=g^2-f^2$ are similar to that of $g^1,f^1,g^1$.  Hence we have
\be \begin{split}H^2(G^2)-H^2(\Sigma^2\bs C^2(o_E,\frac 14 r_E))&\ge\frac14\int_{P_\a^2\cap D(0,\frac 12)\bs D(o_E,\frac 14 r_E)}|\nabla h|^2-Cr_E^2(\mu+\e'\mu+C_0(\mu))\\
&\ge -Cr_E^2(\mu+\e'\mu+C_0(\mu)).\end{split}\ee

Now we still have to estimate the part inside $D(o_E,\frac14 r_E)$. For this purpose we need the following lemma.

\begin{lem}[cf.\cite{2p} Corollary 2.45]\label{bigproj} Suppose $\xi>0$ is such that $\arccos(\xi/2)\le\a_1\le\a_2$, and $P^1,\ P^2$ are two planes with characteristic angles $(\alpha_1,\alpha_2)$. Denote by $p^i$ the orthogonal projection on $P^i,i=1,2$. Then if $E$ is a closed 2- rectifiable set satisfying $p^i(E)\supset B(0,1)\cap P^i$, we have
\be H^2(E)\ge \frac{2\pi}{1+\xi}.\ee
\end{lem}

We apply Lemma 6.29 to the part $E\cap D_\a(o_E,\frac 14 r_E)$, and by Proposition 2.28, we get
\be H^2(E\cap D_\a(o_E,\frac 14 r_E))\ge 2\pi (\frac14 r_E)^2\frac{1}{1+2\cos\theta_1'}.\ee

On the other hand, notice that Lip $f^1<C_0(\mu)$ and Lip $f^2<C_0(\mu)$, we have
\be \begin{split}
H^2(\Sigma^i\cap D_\a(o_E,\frac 14 r_E))&=\int_{P_\a^i\cap D_\a(o_E,\frac 14 r_E)}\sqrt{1+S(f)}\\
&\le\int_{P_\a^i\cap D_\a(o_E,\frac 14 r_E)}\sqrt{1+C_0(\mu)^2+C_0(\mu)^4}\\
&\le\int_{P_\a^i\cap D_\a(o_E,\frac 14 r_E)}1+\frac{C_0(\mu)^2+C_0(\mu)^4}{2}\\
&=\pi (\frac14r_E)^2(1+\frac{C_0(\mu)^2+C_0(\mu)^4}{2}),
\end{split}\ee
therefore
\be H^2(\Sigma\cap D_\a(o_E,\frac 14 r_E))\le 2\pi (\frac14r_E)^2(1+\frac{C_0(\mu)^2+C_0(\mu)^4}{2}).\ee

Thus
\be H^2(\Sigma\cap D_\a(o_E,\frac 14 r_E))-H^2(E\cap D_\a(o_E,\frac 14 r_E))\le2\pi (\frac14r_E)^2(\frac{C_0(\mu)^2+C_0(\mu)^4}{2}+2\cos\a_1).\ee

We combine (6.34), (6.28) and (6.27), and get
\be \begin{split}&H^2(E\cap D(0,\frac 12))-H^2(\Sigma)\\
=&\sum_{i=1,2}[H^2(G^i-H^2(\Sigma^1\bs C^1(o_E,\frac14 r_E))]+[H^2(E\cap D_\a(o_E,r_E))-H^2(\Sigma\cap D_\a(o_E,r_E)]\\
\ge& C_5(\delta')r_E^2-Cr_E^2(\mu+\e'\mu+C_0(\mu)) -Cr_E^2(\mu+\e'\mu+C_0(\mu))\\
&-2\pi (\frac14r_E)^2(\frac{C_0(\mu)^2+C_0(\mu)^4}{2}+2\cos\a_1).\end{split}\ee

Notice that $\d'$ is just a constant, depending on $\e'$, where $\e'$ is the parameter for the $\e'$-process, and guarantees the regularity for parts of minimal sets where the $\e'-$process does not stop. Hence it does not depend on $\mu$ or $\a$. Therefore when $\a$ is large enough and $\mu$ is small enough, 
\be H^2(E\cap D_\a(0,\frac 12))-H^2(\Sigma)>0.\ee

Recall that $\Sigma$ contains a deformation of $E$ in $D_\a(0,\frac 12)$, hence (6.36) contradicts the fact that $E$ is minimal.

This contradiction yields that there exists $\theta_1\in]0,\frac\pi 2[$ and $\mu_0>0$ such that for any $\a>\theta_1$, if $E$ is minimal in $B(0,1)$ with $d_{0,1}(E,P_\a)<\e'$, and moreover (6.1) holds, then $E$ contains a point of type $2\P$ in $B(0,\frac{1}{100})$.

Now for guarantee the condition (6.1), we apply Proposition 2.20 again. Set $\lambda=\e_\mu$. Then when $d_{0,1}(E,P_\a)<\lambda$, our $\lambda$-process does not stop before step 1. Then by (2.22), the curves $\gamma^i$ admits Lipschitz constants less than $\mu$. Thus (6.1) holds.

Thus when $d_{0,1}(E,P_\a)\le\lambda$, there exists a point of type $2\P$ in $B(0,\frac{1}{100})$. This completes the proof of Theorem \ref{main}. \qed

\section{Global regularity and local $C^1$ regularity for minimal sets that are near $2\P$ type minimal cones}

In this section we give two useful corollaries of Theorem \ref{main}, concerning global and local regularity for minimal sets that are near $2\P$ type minimal cones.

\begin{thm}Let $\theta_1$ be as in Theorem \ref{main}. Then for any $\a=(\a_1,\a_2)$ with $\a_2\ge \a_1\ge\theta_1$, if $E$ is a 2-dimensional reduced Almgren minimal set in $\R^4$ such that one blow-in limit of $E$ at infinity is $P_\a$ (i.e., there exists a sequence of numbers $r_n\to \infty$, and the sequence of sets $r_n^{-1}(E)$ converge to $P_\a$ under the Hausdorff distance as $n\to\infty$), then $E$ is a $\P_\a$ set.
\end{thm}

\nd By hypothesis, there exists $R>0$ and a $\P_\a$ set $P_\a$ such that $d_{0,R}(E,P_\a)<\lambda$. Then by Theorem \ref{main}, there exists a $2\P$ type point $x\in E$. In particular, the density $\theta(x)$ of $E$ at $x$ is 2, which is equal to the density $\theta_\infty$ of $E$ at infinity. By the monotonicity (cf.\cite{DJT} Proposition 5.16) of the density function $\theta_x(r)=r^{-d}H^d(E\cap B(x,r))$, it has to be constant for $r\in]0,\infty[$. By Theorem 6.2 of \cite{DJT}, $E$ is a minimal cone centered at $x$. As a result, $d_{x,r}(E,P_\a+x)$ is constant for $r\in]0,\infty[$, since $P_\a+x$ is also a cone centered at $x$. But by hypothesis, $d_{x,r}(E,P_\a+x)\to 0$ as $r\to\infty$, hence $d_{x,r}(E,P_\a+x)=0$, which means that $E=P_\a+x$.\qed

\begin{thm}Let $\theta_1$ be as in Theorem \ref{main}. Then there exists a $\e>0$ such that for any $\a=(\a_1,\a_2)$ with $\a_2\ge \a_1\ge\theta_1$, if $E$ is a 2-dimensional reduced Almgren minimal set in $U\subset\R^4$, $B(x,100r)\subset U$, and there is a reduced minimal cone $P_\a+x$ of type $\P_\a$ centered at $x$ such that $d_{x,100r}(E,P_\a)\le\e$, then there exists a minimal cone $P_{\a'}$ of type $2\P$ such that there is a $C^1$ diffeomorphism $\Phi: B(x,2 r)\to \Phi(B(x,2r))$, such that $|\Phi(y)-y|\le 10^{-2}r$ for $y\in B(x,2r)$, and $E\cap B(x,r)=\Phi( P_{\a'})\cap B(x,r)$.
\end{thm}

\nd Let $\lambda$ be the $\lambda$ in Theorem \ref{main}. Let $\e=\min\{\frac{1}{1000}\lambda,\e_1\}$, where $\e_1$ is the one in Corollary 12.25 of \cite{DEpi}. Then by Theorem \ref{main}, $d_{x,r}(E,P_\a)\le 200d_{x,100r}(E,P_\a)\le\frac{1}{5}\lambda$ yields that there exists a point $y\in B(x,\frac{1}{100}r)$ of type $P_{\a'}$ for some angle $\a'$. 

But $P_{\a'}\cap \partial B(0,1)$ is a disjoint union of two circles, and circles verifies the property of full length because of angles, hence by Remark 14.40 of \cite{DEpi}, $P_{\a'}$ is a minimal cone with the full length property because of angles. We apply Theorem 1.15 of \cite{DEpi}, and get the conclusion. \qed
 
\renewcommand\refname{References}
\begin{spacing}{1}
\bibliographystyle{plain}
\bibliography{reference}
\end{spacing}

\end{document}